\documentclass[a4paper,10pt]{article}
\usepackage{mymath}
\usepackage{tikz}
\usetikzlibrary{external}
\usepgfplotslibrary{groupplots}
\pgfplotsset{compat=1.14}

\title{On the adaptive spectral approximation of functions using redundant sets and frames}
\author{Vincent Copp\'e\footnote{Email: \texttt{vincent.coppe@cs.kuleuven.be}. Website: \texttt{https://people.cs.kuleuven.be/\textasciitilde vincent.coppe}.} \and Daan Huybrechs\footnote{Email: \texttt{daan.huybrechs@cs.kuleuven.be}. Website: \texttt{https://people.cs.kuleuven.be/\textasciitilde daan.huybrechs}.}}

\date{KU Leuven \\ Department of Computer Science \\ Celestijnenlaan 200A \\ 3001 Leuven, Belgium \\ \vskip10pt \today}
\usepackage{makecell}
\def\H{\mathrm{H}}
\def\G{\mathrm{G}}
\def\b{\mathbf{b}}
\def\c{\mathbf{c}}
\def\x{\mathbf{x}}
\def\y{\mathbf{y}}
\def\z{\mathbf{z}}
\def\functionset{\Phi=\left\{\phi_k\right\}_{k\in\NN_0}}
\newcommand{\vertiii}[1]{{\left\vert\kern-0.25ex\left\vert\kern-0.25ex\left\vert #1
		\right\vert\kern-0.25ex\right\vert\kern-0.25ex\right\vert}}
\newcommand{\BB}{\Xi}

\DeclareMathOperator*{\argmin}{arg\,min}

\begin{document}

\maketitle

\begin{abstract}
The approximation of smooth functions with a spectral basis typically leads to rapidly decaying coefficients where the rate of decay depends on the smoothness of the function and vice-versa. The optimal number of degrees of freedom in the approximation can be determined with relative ease by truncating the coefficients once a threshold is reached. Recent approximation schemes based on redundant sets and frames extend the applicability of spectral approximations to functions defined on irregular geometries and to certain non-smooth functions. However, due to their inherent redundancy, the expansion coefficients in frame approximations do not necessarily decay even for very smooth functions. In this paper, we highlight this lack of equivalence between smoothness and coefficient decay and we explore approaches to determine an optimal number of degrees of freedom for such redundant approximations.
\end{abstract}

\section{Introduction}

The approximation of a function $f \in \H$ in a Hilbert space $\H$ using a spectral basis $\Phi \triangleq \left\{\phi_k\right\}_{k\in\NN_0}$ for $\H$, such as a family of orthogonal polynomials in $L^2([-1,1])$, leads to an infinite expansion that can be truncated after finitely many terms $N$:
\begin{equation}\label{eq:expansion}
f(x) \approx f_N(x) = \sum_{k=1}^N c_k\phi_k(x).
\end{equation}
For spectral bases, the coefficients $c_k$ typically decay rapidly if $f$ is smooth. The optimal number of degrees of freedom $N$ can then be determined with relative ease by truncating the coefficients once the coefficient size reaches a given threshold.

A popular spectral basis for function approximation is the sequence of Chebyshev polynomials $\{ T_n \}_{n=0}^\infty$~\cite{mason2003ChebyshevPolynomials}. They are extensively used in software packages such as Chebfun~\cite{driscoll2014fun,trefethen2013approximation} and ApproxFun~\cite{olver2018approxfun}. The optimal coefficient vector length in ChebFun is decided by the method called \texttt{standardChop}\footnote{This is the name of the routine in  Chebfun version 5.3.}. This routine is more involved than a simple truncation based on coefficient size, but it is similar in spirit~\cite{Aurentz2017}. A list of some existing truncation techniques for Chebyshev interpolants is presented in \cite[Chapter 3]{boyd2014solving}.

For any orthonormal basis $\functionset$, the coefficients in \eqref{eq:expansion} are given by $c_k=\langle f,\phi_k\rangle$. \emph{Some} decay of coefficients in orthonormal bases is therefore guaranteed by Bessel's inequality \cite[eqn (2.1.17)]{szeg1939orthogonal} %\cite[Theorem 4.2]{saxe2013beginning},
\begin{equation}\label{eq:besselsinequality}
\sum_{k=1}^N \left|\left\langle f, \phi_k\right\rangle\right|^2 \leq \|f\|_{\H}^2, \quad\forall N\in\NN_0.
\end{equation}
This inequality implies that the series $\sum_{k=1}^\infty \left|\left\langle f, \phi_k\right\rangle\right|^2$ is convergent and bounded by $\|f\|_{\H}^2$. Therefore, the expansion coefficients associated with an orthonormal basis satisfy at least
\[
\lim_{k\rightarrow\infty}c_k = \lim_{k\rightarrow\infty}\left\langle f, \phi_k\right\rangle\ = 0.
\]
Hence, a truncation strategy can in principle be used for expansions in any orthonormal basis.

However, orthonormal bases are not very flexible, and they might not be available for many problems \cite[Ch. 4]{christensen2016}. For example, it may be difficult or even impossible to construct a basis offering spectral approximation accuracy on a geometrically complicated bounded domain $\Omega$. On the other hand, one can easily create a bounding box $\BB$ such that $\Omega\subset\BB$ and create a basis $\Phi$ for $L^2(\BB)$ --- e.g., a tensor product of Fourier series. The basis $\Phi$ for $L^2(\BB)$ is not a basis for $L^2(\Omega)$. Instead, the restriction of an orthonormal basis to a subdomain yields a \emph{frame}~\cite{Adcock2019,christensen2016}. We will refer to such frames obtained from Fourier series on a bounding box as \emph{Fourier extension} frames. Other types of frames that are easily created include the augmentation of an orthonormal basis by a finite number of additional functions which are bounded, and the concatenation of two or more orthonormal bases. This enables the spectral approximation of wide classes of functions including certain non-smooth functions.

A frame is more general and more flexible than a basis. Precise definitions are given in \S\ref{s:preliminaries}. Importantly, like a basis, a frame $ \{ \phi_k \}_{k \in \mathbb{N}_0}$ for a Hilbert space $\H$ is complete, such that any function $f \in \H$ can be written as an expansion $f = \sum_{k \in \mathbb{N}_0} c_k \phi_k$. However, unlike a basis, a frame may be redundant, such that the expansion coefficients $\{c_k\}_{k\in\mathbb{N}_0}$ are not necessarily unique.

Many types of frames that have been studied in signal processing~\cite{Kovacevic2007,Lovisolo2014} or in the field of wavelets~\cite{Daubechies1992} are associated with a so-called \emph{canonical dual frame} $\{ \tilde{\phi}_k \}_{k \in \mathbb{N}_0}$, which has the property that
\begin{equation}\label{eq:dualexpansion} 
 f = \sum_{k \in \mathbb{N}_0} \langle f, \tilde{\phi}_k \rangle \, \phi_k,\quad f\in \H.
\end{equation}
This expansion satisfies an analogue of \eqref{eq:besselsinequality}. Thus, in many cases of interest, it can be truncated after finitely many terms, and result in an approximation to $f$. However, this is not always the preferred approach. Sometimes, the dual frame is not known or can not easily be computed. Even if it  is known, the expansion \eqref{eq:dualexpansion} may not converge rapidly after truncation with increasing $N$~\cite{Adcock2019}.

Function approximation using the truncation of an infinite frame, i.e., using $\Phi_N \triangleq \{ \phi_k \}_{k=1}^N$, was analysed in~\cite{Adcock2019,adcock2017frames}. The approximation coefficients are computed using a regularised singular value decomposition of an associated linear system, which is ill-conditioned if $\Phi_N$ is redundant. Loosely speaking, the ill-conditioning reflects the fact that the coefficients are not unique. Yet, it was shown that highly accuracy and numerically stable computation of an approximation
\[
 f \approx f_N = \sum_{k=1}^N c_{N,k} \phi_k
\]
can be achieved. Fast algorithms are known for univariate and multivariate Fourier extensions~\cite{lyon2011fastcontinuation,matthysen2016fast,matthysen2018function}, and for extension approximations based on B-splines~\cite{coppe2019splines} and wavelets~\cite{waveletextension}. However, as reflected in our notation, the coefficients $\c_N$ depend on $N$: all coefficients may change in nearly arbitrary ways if $N$ is increased, and in practice they usually do. In fact, since the algorithms suggested in~\cite{lyon2011fastcontinuation,matthysen2016fast,matthysen2018function,coppe2019splines,waveletextension} involve routines from randomised linear algebra, the coefficients actually exhibit a fair degree of randomness. Also, perhaps surprisingly, in spite of near-geometric convergence of Fourier extension approximation $f_N$ to $f$ for analytic $f$ and increasing $N$~\cite{huybrechs2010fourier,adcock2014numerical,webb2019pointwise,geronimo2020fourier}, the coefficients do not necessarily decay. 

The computation of spectral approximations using redundant sets with coefficient decay has not yet received a lot of attention. Lyon proposes the solution of a weighted least squares problem for univariate Fourier extension in~\cite{lyon2012smoothing}, and adapts the fast method of~\cite{lyon2011fastcontinuation} for an efficient implementation. Two algorithms were proposed by Gruberger and Levin in~\cite{Gruberger2020}, again for Fourier extension, that achieve coefficient decay. The first algorithm uses weighted least squares, similar to~\cite{lyon2012smoothing}. It is stated in \cite[Lemma 2.3]{Gruberger2020} that in case a good Fourier extension exists, i.e., an extension with rapidly decaying Fourier coefficients, then it will be found numerically. The second method is based on Hermite interpolation (interpolation of function values and derivatives at the endpoints). A similar approach using Hermite interpolation was explored earlier in~\cite{vanbuggenhout2017fourierextension}. The algorithms using weights impose a decay rate a priori, and therefore assume prior knowledge of the smoothness of the function. The algorithms based on Hermite interpolation yield good results for univariate approximations, but are difficult to extend to higher-dimensional domains with general shape. Furthermore, their stability depends on the stability of the Hermite interpolation problem.

Compared to these references, the aim in this paper is the truncation at an optimal number of degrees of freedom. We only assume an algorithm to compute approximation coefficients $\c_N$, we do not assume or enforce that these coefficients decay. (In the process we do suggest one way to do so, without assuming a priori knowledge of $f$.) We aim for approximations using the truncation of an infinite frame or, more generally, an infinite dictionary\footnote{A \emph{dictionary} is a general term in approximation theory to describe a collection of functions that has no discernible structure or properties like a basis or frame, most often used when that collection is overcomplete~\cite{chen1999pursuit}.} $\Phi \triangleq \{ \phi_k \}_{k \in \NN_0}$. In view of the above discussion, truncation based on coefficient decay is not a suitable strategy. Our criterion is to adaptively determine a value of $N$ such that
\begin{equation}\label{eq:criterion_part1}
\|f-f_N\|_\H \leq \delta \|f\|_\H,
\end{equation}
for a given desired relative accuracy $\delta$. The smallest value of $N$ for which the condition holds is the optimal value.
In addition, for reasons outlined further on but mainly for improved stability, we aim for the coefficients to satisfy a bound
\begin{equation}\label{eq:criterion_part2}
\| \c_N \| \leq \mu \|f\|_\H,
\end{equation}
for a given parameter $\mu > 0$. For expansions in orthonormal bases \eqref{eq:criterion_part1} typically implies \eqref{eq:criterion_part2}, but for approximation in the presence of redundancy that is not necessarily the case.

In our implementation, \eqref{eq:criterion_part2} is enforced \emph{implicitly} rather than explicitly since it is not trivial to choose an appropriate value for $\mu$ without intricate knowledge of (frame properties of) the dictionary $\Phi$, and we would like to avoid such expert knowledge. Furthermore, the continuous norms in \eqref{eq:criterion_part2} are replaced by discrete approximations, and a considerable part of the paper is devoted to a justification of this approximation. Fortunately, compared to the many subtleties of the analysis, the algorithms we arrive at are relatively straightforward to use and implement.

We focus on optimal truncations of an infinite frame after $N$ terms. This is in contrast to existing literature on adaptive computations using sparse approximations or best $N$-term approximations in which frames are also sometimes used (see, e.g., \cite{Cohen2000,Cohen2002,Stevenson2003,Dahlke2007}). The best $N$-term approximation $u_N$ to $u$ is based on an optimal subset of $N$ coefficients out of an infinite set. Such adaptive approximations to the solutions of an operator equation $Lu=g$ were explored in~\cite{Cohen2000,Cohen2002} using Riesz bases of wavelet type and extended in~\cite{Stevenson2003} to wavelet frames and in~\cite{Dahlke2007} to Gelfand frames. Efficient algorithms in this context are possible owing to the compression properties of wavelets. In this paper, we fix the subset $1,2,\ldots,N$ of the first $N$ coefficients. This simpler approach may not yield the optimal results of best $N$-term approximations, but is closer in spirit to the truncation strategies of spectral methods. The flexibility of frames and dictionaries enables the extension of spectral methods to much wider classes of functions.

In \S\ref{s:preliminaries}, we briefly recall frames and frame approximations, and we highlight the particular case of Fourier extension (or \emph{Fourier continuation}). This approximation scheme was first proposed in~\cite{boyd2002comparison,bruno2007accurate}. In \S\ref{s:approximationcriterion}, we illustrate the problems associated with truncation based on coefficient decay. We discuss an inefficient way to calculate a frame approximation with decaying coefficients such that truncation can be used, and explore an alternative truncation strategy based on the residual of a least squares problem. This strategy forms the basis of the algorithms in \S\ref{s:algorithms} that determine the optimal approximation of length $N$ adaptively. The first algorithm guarantees optimal $N$ but is costly, the second is much more efficient but only approximates the optimal $N$.

\section{Preliminaries\label{s:preliminaries}}

We recall the definition of frames for a Hilbert space $\H$, introduce examples and quote the relevant theory of function approximation with truncated frames~\cite{Adcock2019,adcock2017frames}. This theory motivates the stopping criterion \eqref{eq:criterion_part1}--\eqref{eq:criterion_part2} for the adaptive frame approximation, as well as the use of a regularised solver for the ill-conditioned linear system that is used to construct the approximation.

We hasten to add that the algorithms of the paper apply even when $\Phi \triangleq \{ \phi_k \}_{k \in \NN_0}$ is not a frame. However, in that case, one can not guarantee the existence of stable approximations for all functions in $\H$. In particular, we can not always guarantee \eqref{eq:criterion_part2}.

\subsection{Frames}

Frames are a generalisation of orthogonal and Riesz bases. A frame for a separable Hilbert space $\H$ is defined as a sequence $\functionset$ that satisfies the so-called \emph{frame condition}
\begin{equation}
A\|f\|_\H^2 \leq \sum_{k\in\NN_0} \left|\left\langle f, \phi_k \right\rangle\right|^2\leq B\|f\|_\H^2,\quad\forall f\in\H,\label{eq:framecondition}
\end{equation}
with constants $0 < A, B < \infty$. A finite upper frame bound $B$ implies boundedness of the associated \emph{Gram operator}. A positive lower frame bound $A > 0$ implies that, like a basis, the frame is complete in $\H$. The frame condition ensures that the samples $\langle f, \phi_k  \rangle$ carry sufficient information to reconstruct the underlying function $f$ (the associated \emph{frame operator} is bounded and boundedly invertible). See~\cite{christensen2016} for details and the preliminary section of~\cite{Adcock2019} for a concise overview.

Unlike a basis, the frame condition~\eqref{eq:framecondition} allows the frame elements $\phi_k$ to be linearly dependent, or nearly linearly dependent. This can be seen with a simple example. Consider the concatenation of two different orthonormal bases $\Phi=\{\phi_k\}_{k\in\NN_0}$ and $\Psi = \{\psi_k\}_{k\in\NN_0}$ for the same space $\H$. This concatenated set is clearly complete, because both $\Phi$ and $\Psi$ are complete. Applying the Riesz identity twice, i.e., $\Vert f \Vert_\H^2 = \sum_{k\in\NN_0} \left|\left\langle f, \phi_k \right\rangle\right|^2 = \sum_{k\in\NN_0} \left|\left\langle f, \psi_k \right\rangle\right|^2$, we see that the set satisfies the frame condition \eqref{eq:framecondition} with $A=B=2$. To see that there are an infinite number of representations for every $f\in\H$, it is enough to note that the zero function can be represented in infinitely many ways. For all $g\in\H$, we can find the basis coefficients $\{a_k(g)\}_{k\in\NN_0}$ and $\{b_k(g)\}_{k\in\NN_0}$ for the bases $\Phi$ and $\Psi$ respectively. Then, for each $g\in\H$, the concatenation of coefficients $\{a_k(g)\}_{k\in\NN_0}\cup \{-b_k(g)\}_{k\in\NN_0}$ consists of exact coefficients in the frame $\Phi\cup\Psi$ of the zero function.

\subsubsection{Example: Fourier extensions}

\begin{figure}[tb]
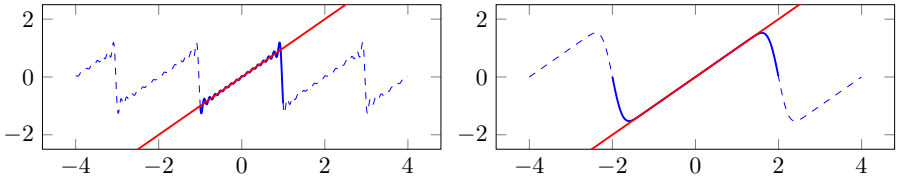

	\centering
	\resizebox{.8\columnwidth}{!}{\input{img/1Dexample_Gibbs.tikz}\input{img/1Dexample.tikz}}
	\caption{Approximation of $f(x)=x$ on $[-1,1]$ by a classical Fourier series on $[-1,1]$ (left panel) and by a Fourier extension frame on $[-2,2]$ (right panel). The frame approximation does not suffer from the Gibbs phenomenon.}\label{fig:1dexample}
\end{figure}

A Fourier extension frame is a type of frame suitable for spectral function approximation on domains with arbitrarily complex geometries. In order to approximate a function on a bounded domain $\Omega$, this domain is embedded into a bounding box~$\BB$, $\Omega\subset\BB$. A tensor product Fourier series on $\BB$ is an orthogonal basis for $L^2(\BB)$, which is used to approximate functions in $L^2(\Omega)$. Implicitly, this corresponds to extending a function from $\Omega$ to a periodic function on $\BB$. One can think of the redundancy in this frame as corresponding to the many different ways in which a smooth function can be extended to another smooth function on a larger domain.

A 1-D Fourier extension example is shown in \Cref{fig:1dexample}. The function $f(x) = x$ is smooth but non-periodic on $[-1,1]$ and famously gives rise to the Gibbs phenomenon when approximated by a Fourier series on $[-1,1]$~\cite{gibbs1898fourier}. In contrast, when using Fourier series on $\BB = [-2,2]$, periodicity on $\Omega = [-1,1]$ is not required and extensions exist that exhibit no Gibbs phenomenon and offer rapid numerical convergence. One such extension is shown in the right panel of \Cref{fig:1dexample}.

\subsection{Computing frame and dictionary approximations}\label{ss:frameapproximation}

We summarise the methodology of~\cite{Adcock2019}.
Suppose we want to approximate $f \in \H$ using the truncation $\Phi_N=\left\{\phi_n\right\}_{n=1}^{N}$ of an infinite frame for $\H$. We may look for the approximant $f_N \in \H_N \triangleq \myspan\{\Phi_N\}$ that is closest to $f$ in the $\H$-norm:
 \[
 f_N = \argmin_{g\in \H_N}\|f-g\|_{\H}.
 \]
This is equivalent to finding the orthogonal projection $\mathcal P_Nf$ of $f$ onto $\H_N$. An orthogonal projection in general (unless $\Phi_N$ is an orthogonal basis) requires the solution of a linear system
\begin{equation}\label{eq:system}
	G_N \c = \b, \quad \b=\left\{\left\langle f, \phi_n \right\rangle\right\}_{n=1}^N,
\end{equation}
where $G_N$ is the truncated $N\times N$ Gram-matrix
\[
 G_N = \left\{\left\langle \phi_i,\phi_j \right\rangle\right\}_{i,j=1}^{N}\in\CC^{N\times N}.
\]
For an orthogonal basis, the Gram matrix is diagonal.

Because a frame is redundant, the conditioning of the linear system \eqref{eq:system} can be arbitrarily bad. Despite this potential ill-conditioning, regularisation may lead to accurate and stable approximations. The regularisation proposed in~\cite{Adcock2019} consists of truncating the singular values of the singular value decomposition of $G_N$ below a threshold $\epsilon$ in solving the system \eqref{eq:system}. Denoting $\mathcal{P}_N^\epsilon$ as the solution obtained from the regularised projection, a generic error bound can be shown to hold for all $f \in \H$ \cite[Theorem 13]{Adcock2019}:
\begin{equation}\label{eq:fna1bound}
	\|f-\mathcal{P}_N^\epsilon f\|_{\H}\leq \inf_{\z\in\CC^N} \left\{ \|f - \mathcal{T}_N \z\|_{\H}+\sqrt{\epsilon} \|\z\| \right\},
\end{equation}
where $\mathcal{T}_N \z = \sum_{n=1}^Nz_n\phi_n\in\H_N$.
A suitable solution is found numerically if some coefficient vector $\z$ exists which gives a good approximation to $f$ and such that the discrete $\ell_2$ norm $\|\z\|$ is small.

The error bound \eqref{eq:fna1bound} holds for all sets $\Phi_N$ of bounded functions in $\H$. In particular, the set $\Phi_N$ does not need to arise from the truncation of a frame, the infinite set $\Phi$ may be any dictionary. The frame condition ensures the existence of vectors that make the right hand side of \eqref{eq:fna1bound} small for any $f \in \H$ by taking $N$ sufficiently large. To be precise, the infinite dual frame expansion \eqref{eq:dualexpansion}, which --- although it might not be easily computed in practice --- is guaranteed to exist, has coefficients with norm bounded by $\frac{1}{\sqrt{A}}\Vert f \Vert_{\H}$, where $A$ is the lower frame bound in \eqref{eq:framecondition}. Thus, if $\Phi_N$ arises from a frame, an adaptive approximation is guaranteed to eventually succeed. If $\Phi_N$ is not a frame, then the adaptive approximation may succeed for some $f\in\H$ and fail for some others.

\subsection{Discrete approximation and the stable sampling rate}\label{ss:discrete_approximation}

The error bound~\eqref{eq:fna1bound}  yields accuracy up to $\sqrt{\epsilon}$ only, and the procedure with the Gram matrix requires the computation of a large number of inner products. This situation is improved by considering generalised sampling in combination with oversampling in the follow-up paper~\cite{adcock2017frames}.

Assume that the data on $f \in \H$ is given by $M$ functionals $\{ l_{m,M}(f) \}_{m=1}^M$ applied to $f$. Two examples of families of functionals include inner products: $l_{m,M}(f)=\langle f, \phi_m\rangle$, and point evaluations: $l_{m,M}(f)=f(x_{m,M})$ for a set of $M$ points $\{x_{m,M}\}_{m=1}^M$. Sampling with inner products and choosing $M=N$ yields the Gram matrix again. Oversampling corresponds to choosing $M > N$.

Numerically stable approximations may again be feasible in spite of redundancy in the truncated frame, albeit with some additional conditions. The expansion coefficients are found by solving the least squares problem
\begin{equation}\label{eq:collocationsystem}
 A \c = \b, \quad \b = \left\{l_{m,M}(f)\right\}_{m=1}^M,\quad A = \left\{ l_{m,M}(\phi_n) \right\}_{m=1,n=1}^{M,N}\in\CC^{M\times N}.
\end{equation}
We regularise using a truncated singular value decomposition as above, i.e.,
\[
 A = U \Sigma V^T \approx U \Sigma_\epsilon V^T =: A_\epsilon,
 \]
where $\Sigma_\epsilon$ is $\Sigma$ with all singular values smaller than $\epsilon$ replaced by $0$. We denote the regularised solution by
\begin{equation}
 P^\epsilon_{M,N}f = A_\epsilon^\dagger \b.
\end{equation}
This yields the error bound \cite[Theorem 1.3]{adcock2017frames}:
\begin{equation}\label{eq:fna2bound}
\|f-\mathcal P^\epsilon_{M,N}f\|_{\H}\leq \inf_{\z\in\CC^N} \left\{ \|f-\mathcal T_N \z\|_{\H}+\kappa_{M,N}^\epsilon\|f-\mathcal T_N \z\|_M+\epsilon\lambda_{M,N}^\epsilon\|\z\| \right\}.
\end{equation}
Here, $\kappa_{M,N}^\epsilon$ and $\lambda_{M,N}^\epsilon$ are constants that depend on the sampling functionals $l_{m,M}$ and on the frame elements $\phi_n$. The $M$-norm $\Vert \cdot \Vert_M^2 \triangleq \sum_{m=1}^m | l_{m,M}(\cdot) |^2$ is a data-dependent norm. Combining its definition with~\eqref{eq:collocationsystem} one can conclude that
\begin{equation}\label{eq:sampling_error}
 \Vert f - \mathcal T_N\z\Vert_M = \Vert A \z- \b \Vert,
\end{equation}
i.e., the $M$-norm appearing in the right hand side of \eqref{eq:fna2bound}, but absent from the earlier bound \eqref{eq:fna1bound}, is nothing but the discrete $\ell_2$ norm of the residual vector of the linear system.

The analysis of the $\kappa_{M,N}^\epsilon$ and $\lambda_{M,N}^\epsilon$ constants that appear in \eqref{eq:fna2bound} is highly involved and frame-specific. However, if the sampling functionals are sufficiently `rich' for increasing $M$, a concept that is defined below, it is guaranteed that both constants are bounded for $M$ sufficiently large relative to $N$ \cite[Proposition 4.6]{adcock2017frames}. Furthermore, it is shown that there exists a \emph{stable sampling rate} $M = \Theta^\epsilon(N,\theta)$, with $1 < \theta < \infty$, such that for $M \geq \Theta^\epsilon(N,\theta)$ both constants satisfy the bound
\begin{equation}\label{eq:bound}
 \kappa_{M,N}^\epsilon,\lambda_{M,N}^\epsilon \leq \frac{\theta}{A'}, \qquad 1 < \theta < \infty.
\end{equation}
The constant $A'$ appearing in the denominator is determined by the family of sampling functionals. It arises from the `richness' condition \cite[(1.7)]{adcock2017frames} with $0 < A' \leq B' < \infty$,
\begin{equation}\label{eq:richness}
 A' \Vert f \Vert_{\H}^2 \leq \liminf_{M \to \infty} \Vert f \Vert_M^2 \leq \limsup_{M \to \infty} \Vert f \Vert_M^2  \leq B' \Vert f \Vert_{\H}^2, \quad \forall f \in \G.
\end{equation}
Here, $\G \subset \H$ is a subspace of $\H$ on which the sampling functionals are well-defined. In the main examples of this paper $\H = L^2(\Omega)$ and $\G = L^\infty(\Omega)$, since the latter space allows point evaluations (and hence discrete sampling) but the former does not. On the other hand, $\H$ is a Hilbert space but $\G$ is not.

Finally, we note that the coefficients $\mathbf{c}^\epsilon$ of the approximation $\mathcal P^\epsilon_{M,N} f = \sum_{k=1}^N c_k^\epsilon \phi_k$ can also be bounded \cite[Theorem 4.5]{adcock2017frames}
\begin{equation}\label{eq:coefficient_bound}
 \Vert \mathbf{c}^\epsilon \Vert \leq \inf_{\z \in \mathbb{C}^N} \left\{ \frac{1}{\epsilon} \Vert f - \mathcal T_N \z \Vert_{M} + \Vert \z \Vert \right\}.
\end{equation}
This bound implies that, unless $f$ can be well approximated in $\Phi_N$ such that the first term above is small, the coefficients found by solving \eqref{eq:collocationsystem} may be as large as $\epsilon^{-1}$. As soon as $N$ is sufficiently large to resolve $f$, the first term in the bound may decrease and the coefficient norm settles down. The lower limit is given by the norm of the dual frame expansion  coefficients \eqref{eq:dualexpansion}, that satisfies $\Vert \c \Vert \leq \frac{1}{\sqrt{A}}\Vert f \Vert_{\H}$ where $A$ is the lower frame bound.

\subsubsection{Practical choices}

In practice, we will assume that the stable sampling rate is known. For all examples in this paper, we simply choose a linear oversampling rate
\begin{equation}\label{eq:linearoversampling}
 M = \gamma N, \qquad \gamma > 1,
\end{equation}
for a value of $\gamma$ that is deemed sufficiently large based on experiments for a given frame (typically $\gamma=2$). Furthermore, all our examples are based on the function space $\H = L^2(\Omega)$ for some domain $\Omega$, and as sampling functionals we choose weighted point evaluations of the form
\[
 l_{m,M}(f) = w_{m,M} f(x_{m,M}), \qquad f \in G = L^\infty(\Omega).
\]
The points belong to the domain $\Omega$ at hand, and the weights $w_{m,M} > 0$ are chosen such that
\begin{equation}\label{eq:riemann_sum}
 \Vert f \Vert_M^2 = \sum_{m=1}^M w_{m,M}^2 |f(x_m,M)|^2 \to \int_\Omega |f(x)|^2 {\rm d}x = \Vert f \Vert_{L^2(\Omega)}, \qquad M \gg 1.
\end{equation}
This is realised by forming a Riemann sum for $\Vert f \Vert_{L^2(\Omega)}$ in the limit $M \to \infty$. We stress that the agreement between $\Vert f \Vert_M$ and $\Vert f \Vert_{\H}$ need not be highly accurate.  As it will merely be used to detect convergence later on, it does not play a role in the approximation accuracy. Still, we have conveniently achieved that $A'=B'=1$. When using equispaced points, a suitable choice is
\[
w_{m,M} = \sqrt{\frac{|\Omega|}{M}}.
\]
Importantly, these choices do not restrict the domain shape  $\Omega$. Even for a very irregular geometry, an equispaced grid on a bounding box can be restricted to $\Omega$ and the weights defined above yield a Riemann sum in the limit $M \to \infty$.

Furthermore, in the adaptive scheme we will not assume that the user knows exactly what the frame bounds $A$ and $B$ are, nor what the constants $A'$ and $B'$ are in \eqref{eq:richness}.

\section{Truncation strategies for frame approximations}\label{s:approximationcriterion}

An essential element of any adaptive approach for finding the optimal function set size $N$ is the criterion that gives an indication on whether or not the approximation $f_N$ of size $N$ is sufficiently accurate. We say that a solution is accurate if, given $\delta>0$,
\begin{equation}\label{eq:simplestopcriterion}
\|f-f_N \|_\H\leq \delta \|f\|_\H.
\end{equation}
As mentioned in the introduction, a stopping criterion based on coefficient decay is possible in general when approximating with a spectral basis, but that is no longer the case when approximating with truncated frames. We illustrate this restriction in \S\ref{ss:coefsizeonly} and consider an alternative criterion in \S\ref{ss:residual}--\S\ref{ss:stopcriterion}.

\subsection{Truncation based on coefficient size}\label{ss:coefsizeonly}

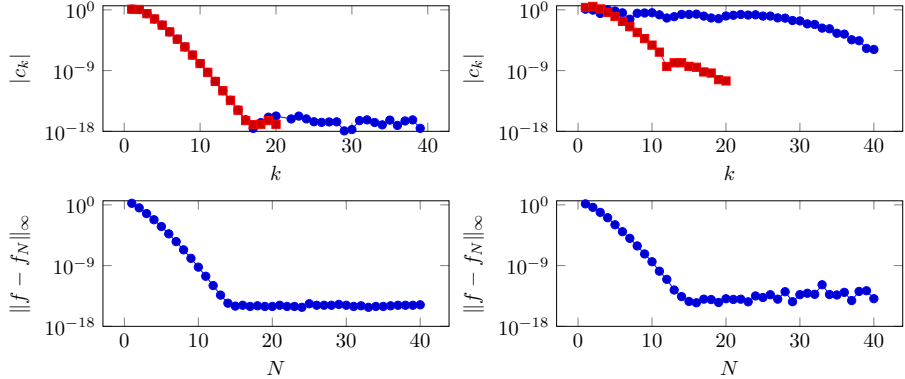
\begin{figure}[t]
	\centering
	\resizebox{.8\columnwidth}{!}{\begin{tikzpicture}
\begin{axis}[width={.5\textwidth}, height={.25\textwidth}, ymode={log}, ymin={1.0e-18}, ymax={4}, ylabel={$|c_k|$}, xlabel={$k$}]
    \addplot+[]
        table[row sep={\\}]
        {
            \\
            1.0  1.2660658777520084  \\
            2.0  1.1303182079849698  \\
            3.0  0.2714953395340765  \\
            4.0  0.04433684984866382  \\
            5.0  0.00547424044209379  \\
            6.0  0.0005429263119140767  \\
            7.0  4.49773229542679e-5  \\
            8.0  3.198436462396476e-6  \\
            9.0  1.992124806862671e-7  \\
            10.0  1.1036771749985444e-8  \\
            11.0  5.505896481558649e-10  \\
            12.0  2.497951729972941e-11  \\
            13.0  1.0391188065635578e-12  \\
            14.0  3.9919556958225997e-14  \\
            15.0  1.4241654844235425e-15  \\
            16.0  5.241080267573207e-17  \\
            17.0  2.7423322473900757e-18  \\
            18.0  1.987260068158088e-17  \\
            19.0  1.2461816042730103e-16  \\
            20.0  1.8391612474776623e-16  \\
            21.0  0.0  \\
            22.0  6.811410945660194e-17  \\
            23.0  1.800988599538586e-16  \\
            24.0  7.004084830718572e-17  \\
            25.0  2.620428421725814e-17  \\
            26.0  2.1219967131617728e-17  \\
            27.0  2.452232313193216e-17  \\
            28.0  2.536849450212524e-17  \\
            29.0  1.1644704110484651e-18  \\
            30.0  1.8491243525182533e-18  \\
            31.0  3.863016244858344e-17  \\
            32.0  4.9460399827559803e-17  \\
            33.0  2.0396331023734894e-17  \\
            34.0  9.156256753229537e-18  \\
            35.0  4.600134292584434e-17  \\
            36.0  7.61109925084824e-18  \\
            37.0  3.434752482434078e-17  \\
            38.0  4.926614671774132e-17  \\
            39.0  2.775557561562891e-18  \\
            40.0  0.0  \\
        }
        ;
    \addplot+[]
        table[row sep={\\}]
        {
            \\
            1.0  1.2660658777520084  \\
            2.0  1.13031820798497  \\
            3.0  0.2714953395340766  \\
            4.0  0.04433684984866383  \\
            5.0  0.005474240442093641  \\
            6.0  0.0005429263119139396  \\
            7.0  4.49773229543168e-5  \\
            8.0  3.198436462387002e-6  \\
            9.0  1.9921248059068895e-7  \\
            10.0  1.1036771793074786e-8  \\
            11.0  5.50589660800999e-10  \\
            12.0  2.4979566179581142e-11  \\
            13.0  1.039135095571416e-12  \\
            14.0  3.9883227844178094e-14  \\
            15.0  1.350059387200675e-15  \\
            16.0  3.9508327165371785e-17  \\
            17.0  1.0061396160665482e-17  \\
            18.0  1.1102230246251566e-17  \\
            19.0  3.885780586188048e-17  \\
            20.0  1.1102230246251566e-17  \\
        }
        ;
\end{axis}
\end{tikzpicture}\begin{tikzpicture}
\begin{axis}[width={.5\textwidth}, height={.25\textwidth}, ymode={log}, ymin={1.0e-18}, ymax={4}, ylabel={$|c_k|$}, xlabel={$k$}]
    \addplot+[]
        table[row sep={\\}]
        {
            \\
            1.0  1.249931361193708  \\
            2.0  0.9129454931886687  \\
            3.0  0.28805022857233026  \\
            4.0  1.0251189204047042  \\
            5.0  0.6666345585655098  \\
            6.0  0.35768589850570015  \\
            7.0  0.0413089190939344  \\
            8.0  0.3078598635561769  \\
            9.0  0.34111445868162116  \\
            10.0  0.37040002958771845  \\
            11.0  0.18274196211721488  \\
            12.0  0.06156630793050799  \\
            13.0  0.09775632859771215  \\
            14.0  0.2112240186207126  \\
            15.0  0.2127639591119346  \\
            16.0  0.23274373081387664  \\
            17.0  0.11938474695845187  \\
            18.0  0.06253707097028313  \\
            19.0  0.0489244491372794  \\
            20.0  0.12234367821328183  \\
            21.0  0.15745284658584804  \\
            22.0  0.2058344827700086  \\
            23.0  0.16817074617410907  \\
            24.0  0.18569320695649724  \\
            25.0  0.11986697066721984  \\
            26.0  0.11982774938010836  \\
            27.0  0.06387488395340171  \\
            28.0  0.05930836574240432  \\
            29.0  0.026253711063053173  \\
            30.0  0.022936224866217353  \\
            31.0  0.008317642139136591  \\
            32.0  0.00688988414459824  \\
            33.0  0.0019835825202159218  \\
            34.0  0.0015657301482701863  \\
            35.0  0.0003377628860044272  \\
            36.0  0.00025493698058663076  \\
            37.0  3.6821904857211964e-5  \\
            38.0  2.6640751237940676e-5  \\
            39.0  1.9406780408624124e-6  \\
            40.0  1.3482631089544374e-6  \\
        }
        ;
    \addplot+[]
        table[row sep={\\}]
        {
            \\
            1.0  2.2795849698373307  \\
            2.0  3.1812728919485096  \\
            3.0  1.3778962882970496  \\
            4.0  0.42547923257674847  \\
            5.0  0.10145667908865322  \\
            6.0  0.01965087752380943  \\
            7.0  0.0032000579138500977  \\
            8.0  0.00044899860339300173  \\
            9.0  5.525130968748845e-5  \\
            10.0  5.955551719086785e-6  \\
            11.0  5.434113200466019e-7  \\
            12.0  4.062458669219165e-9  \\
            13.0  1.4270670069849052e-8  \\
            14.0  1.4384997039991511e-8  \\
            15.0  4.239208573699223e-9  \\
            16.0  3.125421190038913e-9  \\
            17.0  6.270148081411189e-10  \\
            18.0  4.305195223160684e-10  \\
            19.0  4.502414542898468e-11  \\
            20.0  2.8973212653043746e-11  \\
        }
        ;
\end{axis}
\end{tikzpicture}}%

	\resizebox{.8\columnwidth}{!}{\begin{tikzpicture}
\begin{axis}[width={.5\textwidth}, height={.25\textwidth}, ymode={log}, ymin={1.0e-18}, ymax={4}, ylabel={$\|f-f_N\|_{\infty}$}, xlabel={$N$}]
    \addplot+[]
        table[row sep={\\}]
        {
            \\
            1.0  1.7033688723475995  \\
            2.0  0.3573440959655372  \\
            3.0  0.05117239453531708  \\
            4.0  0.006240603862831762  \\
            5.0  0.0006176642478790129  \\
            6.0  5.065822359506811e-5  \\
            7.0  3.5695503153299057e-6  \\
            8.0  2.1563215213404874e-7  \\
            9.0  1.1938530164457006e-8  \\
            10.0  5.926730217709064e-10  \\
            11.0  2.6403768060845323e-11  \\
            12.0  1.1075584893660562e-12  \\
            13.0  4.263256414560601e-14  \\
            14.0  2.6645352591003757e-15  \\
            15.0  9.992007221626409e-16  \\
            16.0  1.3322676295501878e-15  \\
            17.0  8.881784197001252e-16  \\
            18.0  1.1102230246251565e-15  \\
            19.0  8.881784197001252e-16  \\
            20.0  8.881784197001252e-16  \\
            21.0  1.3322676295501878e-15  \\
            22.0  8.881784197001252e-16  \\
            23.0  8.881784197001252e-16  \\
            24.0  6.661338147750939e-16  \\
            25.0  2.220446049250313e-15  \\
            26.0  1.3322676295501878e-15  \\
            27.0  1.3322676295501878e-15  \\
            28.0  1.3322676295501878e-15  \\
            29.0  1.7763568394002505e-15  \\
            30.0  1.3322676295501878e-15  \\
            31.0  8.881784197001252e-16  \\
            32.0  1.1102230246251565e-15  \\
            33.0  6.661338147750939e-16  \\
            34.0  8.881784197001252e-16  \\
            35.0  8.881784197001252e-16  \\
            36.0  1.1102230246251565e-15  \\
            37.0  1.3322676295501878e-15  \\
            38.0  1.3322676295501878e-15  \\
            39.0  1.3322676295501878e-15  \\
            40.0  1.5543122344752192e-15  \\
        }
        ;
\end{axis}
\end{tikzpicture}\begin{tikzpicture}
\begin{axis}[width={.5\textwidth}, height={.25\textwidth}, ymode={log}, ymin={1.0e-18}, ymax={4}, ylabel={$\|f-f_N\|_{\infty}$}, xlabel={$N$}]
    \addplot+[]
        table[row sep={\\}]
        {
            \\
            1.0  1.4017157334070058  \\
            2.0  0.42667694069908624  \\
            3.0  0.07728231206734515  \\
            4.0  0.012469373559859331  \\
            5.0  0.0011173119156064892  \\
            6.0  0.00010647009307263033  \\
            7.0  1.103199896457241e-5  \\
            8.0  8.268572679881458e-7  \\
            9.0  5.3087937601503654e-8  \\
            10.0  3.761188249740144e-9  \\
            11.0  1.4864059982855338e-10  \\
            12.0  9.203693362991316e-12  \\
            13.0  2.4175106361212784e-13  \\
            14.0  2.5701663020072374e-14  \\
            15.0  4.884981308350689e-15  \\
            16.0  3.1086244689504383e-15  \\
            17.0  1.071365218763276e-14  \\
            18.0  9.325873406851315e-15  \\
            19.0  3.4416913763379853e-15  \\
            20.0  1.2212453270876722e-14  \\
            21.0  9.769962616701378e-15  \\
            22.0  1.0547118733938987e-14  \\
            23.0  4.274358644806853e-15  \\
            24.0  3.197442310920451e-14  \\
            25.0  1.865174681370263e-14  \\
            26.0  4.829470157119431e-14  \\
            27.0  1.0658141036401503e-14  \\
            28.0  1.3500311979441904e-13  \\
            29.0  4.440892098500626e-15  \\
            30.0  4.973799150320701e-14  \\
            31.0  7.926992395823618e-14  \\
            32.0  5.545564008002657e-14  \\
            33.0  1.4553913629811177e-12  \\
            34.0  5.1958437552457326e-14  \\
            35.0  4.440892098500626e-14  \\
            36.0  1.092459456231154e-13  \\
            37.0  6.661338147750939e-15  \\
            38.0  1.4876988529977098e-13  \\
            39.0  1.8784973576657649e-13  \\
            40.0  1.2545520178264269e-14  \\
        }
        ;
\end{axis}
\end{tikzpicture}}%
	\caption{\label{fig:illustrationunexpectedproblem} Top: The coefficient sizes of the approximation of $f(x)=e^x$ on $[-1,1]$ with Chebyshev polynomials on $[-1,1]$  (left), and with Chebyshev polynomials on $[-2,2]$ (right) with 20 and 40 degrees of freedom in red squares and blue dots respectively. Bottom: The uniform approximation error for increasing approximation size.}
\end{figure}

The decreasing coefficient size in the approximation of a smooth function using an orthogonal basis is shown in the top-left panel of Figure~\ref{fig:illustrationunexpectedproblem}. The panel shows the size of the expansion coefficients for $e^x$ using Chebyshev polynomials on $[-1,1]$. For spectral approximation schemes, the coefficients decay at a rate that increases with increasing smoothness of the function: since $e^x$ is entire, its expansion coefficients decay super exponentially~\cite{boyd2001spectral}.

This property is not seen when approximating with frames. In the right panel of \Cref{fig:illustrationunexpectedproblem}, $e^x$ is approximated on $[-1,1]$ with Chebyshev polynomials that were scaled to $[-2,2]$. The approximation is computed by solving \eqref{eq:collocationsystem} with regularised least squares, using $\epsilon=10^{-14}$. This could be called a Chebyshev extension problem: $f(x)=e^x$ is implicitly extended from $[-1,1]$ to a polynomial of degree $N-1$ on $[-2,2]$. Here, some decay of the coefficients is observed for the two approximations shown, corresponding to $N=20$ and $N=40$. However, the approximation with larger $N$ actually has larger coefficients and slower decay. Thus, the decay rate of the coefficients does not offer clear information on the optimal value of $N$.

Furthermore, the size of the coefficients does not even correspond closely to the approximation error. The approximation errors for the two experiments are shown in the bottom row of the figure. In both cases, machine precision accuracy is roughly reached, and the optimal values of $N$ are quite close for the basis approximation (left) and the frame approximation (right). However, the size of the error in the Chebyshev extension approximation problem is much smaller than the size of its expansion coefficients. All coefficients shown in the top-right panel, including the late ones, are orders of magnitude larger than machine precision.

One may be led to believe that the lack of coefficient decay is an artefact of the formulation of the problem as a least squares approximation in \eqref{eq:collocationsystem}. Indeed, the current experiment shows absence of coefficient decay of only one particular solution to the linear system. In view of the redundancy of the frame, there might be other representations of the same function with rapidly decaying coefficients. No doubt this is the case. Thus, the fact that rapidly decaying coefficients are not recovered by the regularised least squares solver means that such coefficient vectors must yield a larger right hand side in the bound \eqref{eq:fna2bound}. However, the question remains how to compute those approximations. The ill-conditioning of the system matrix (in the current example the matrix is singular to working precision once $N > 10$) does not prevent computation of vectors that correspond to a small residual (backward error), but it does prevent computation of any specific solution vector to high accuracy (forward error). In \S\ref{s:adaptivedecay}, we will describe one way to compute approximations with rapidly decaying coefficients, albeit a costly one.

Here, we illustrate the inherent difficulty with a modification of the linear system, in which we penalise the late coefficients. Thus, rather than solving $A \c^\epsilon= \b$ for the coefficient vector $\c^\epsilon$, we solve a weighted least squares problem
\begin{align}
	AD\y^\epsilon &= \b\label{eq:weightedsystem1}\\
	\c^\epsilon&=D\y^\epsilon\nonumber
\end{align}
using a truncated SVD, with threshold $\epsilon$, in the first step. Matrix $D$ is a diagonal matrix with diagonal entries chosen to be, say, $d_{nn} = n^{-\alpha}$, with $\alpha > 0$. The original least squares problem yields, among all possible vectors with small residual, the vector with smallest norm $\Vert \c^\epsilon \Vert$. The weighted least squares problem yields a vector with small norm $\Vert \y^\epsilon \Vert$, hence $\c^\epsilon=D\y^\epsilon$ may be expected to exhibit $c_k \sim k^{-\alpha}$ decay.  This approach is similar to the method of~\cite{lyon2012smoothing} and to the first algorithm in \cite{Gruberger2020}. One difference is that, in \cite{Gruberger2020}, the ill-conditioned system is solved using iterative refinement. A disadvantage of this weighted scheme is that one has to decide a priori on the decay rate $\alpha$, which assumes knowledge of the smoothness properties of $f$.

The weighted least squares formulation is illustrated in Figure~\ref{fig:smoothing} (blue line) for a Chebyshev extension approximation. The diagonal matrix has entries that decay algebraically with $k$, the degree of the corresponding Chebyshev polynomial. The smoothed Chebyshev extension approximation does have coefficients that decay quicker than in the original non-weighted formulation of the problem. However, the two issues identified above remain: approximations with larger $N$ have larger coefficients, giving no indication of the optimal value of $N$, and the size of the coefficients remains significantly larger than the actual approximation error on $[-1,1]$.

\subsection{Optimal coefficient decay}\label{s:adaptivedecay}

Optimal coefficient decay without a priori knowledge of the smoothness of $f$ can also be achieved with a straightforward, albeit computationally expensive, modification to the weighted least squares problem. It is achieved by choosing the weight of coefficient $c_{i+1}$ proportional to the best approximation error using $i$ degrees of freedom.

Assume that we have solved the frame approximation problem with $i$ degrees of freedom, i.e., we have solved $A_i\c^\epsilon_i=\b_i$ with a residual norm $r_i \triangleq \|A_i\c^\epsilon_i-\b_i\|$. This means that the approximation error $f-f_i$ is on the order of $r_i$, at least pointwise. 
If $\Phi$ is normalised, it is reasonable to assume that we can approximate the tail $f-f_i$ with coefficients of size $r_i$ or smaller. Thus, we can use the residual corresponding to $i$ degrees of freedom as a weight for the later degrees of freedom in the weighted least squares problem. Doing so for $i=1,\ldots,N$ leads to Algorithm~\ref{alg:rapiddecay}.
\begin{algorithm}
	\caption{Incrementally weighted least squares}\label{alg:rapiddecay}
	{\bf Input:} $A_1,\dots,A_N=A$, $\b_1,\dots,\b_N=\b$ \\
	{\bf Output:} $\c^\epsilon$ such that $A\c^\epsilon \approx \b$
	\begin{algorithmic}
		\State Solve $A_1 \c^\epsilon_1 = \b_1$
		\State $e_1 \gets \|A_1 \c^\epsilon_1 - \b_1\|$
		\State $D_2=\mbox{diag}(\|\b\|,e_1)$
		\For{$i=2,\dots,N$}
		\State Solve $A_i D_i \y^\epsilon_i = \b_i$
		\State $e_i \gets \|A_i D_i \y^\epsilon_i - \b_i\|$
		\State $D_{i+1}=\mbox{diag}(\|\b\|,e_1, e_2 \dots,e_{i-1},e_{i})$
		\EndFor
		\State $\c^\epsilon\gets D_N\y^\epsilon_N$
	\end{algorithmic}
\end{algorithm}

\begin{figure}[t]
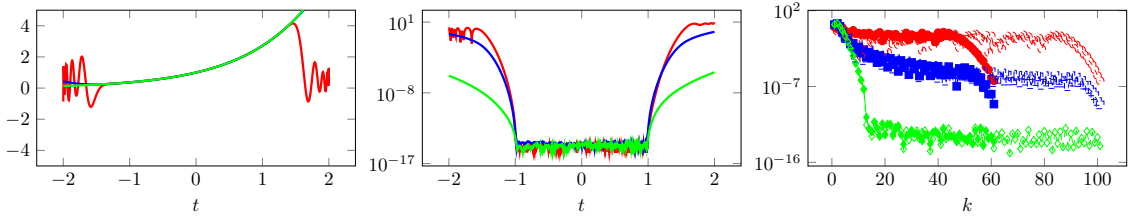

	\centering
	\resizebox{\textwidth}{!}{\input{img/approximation.tikz}\input{img/error.tikz}\input{img/coefficients.tikz}}
	\caption{Approximation of $f(x)=e^x$ in the interval $[-1,1]$ with Chebyshev polynomials on $[-2,2]$. Left: the approximation $\mathcal P^\epsilon_{M,N}f(t)$  with extension, $\epsilon=10^{-14}$. Middle: the pointwise error $|f(t)-\mathcal P^\epsilon_{M,N}(t)|$, $N=61$. Right: the sizes of $\mathbf c_{k,N}$, $N=61$ (full line, filled marks) and $N=101$ (dashed line).
	Red (circle): Approximation, without weight, see \eqref{eq:collocationsystem}. Blue (square): Approximation with weight $D$, $D(k,k) = (10^{-4}+|k|+|k|^2+|k|^3)^{-1}$, see \eqref{eq:weightedsystem1}. Green (diamond): Incrementally weighted least squares, see Algorithm~\ref{alg:rapiddecay}.}\label{fig:smoothing}
\end{figure}

We present no rigorous analysis of this algorithm, but we illustrate its potential with an experiment in Figure~\ref{fig:smoothing} (green line). In this experiment, we obtained coefficients that decay down to the approximation error, at the maximal rate allowed by the smoothness of the function and the accuracy of the regularised solver. The method is generally applicable, but comes at the cost of solving approximation problems with increasing number of degrees of freedom. This could be optimised by using the same weight for a range of degrees of freedom, thereby reducing the number of approximation problems that have to be solved. We do not pursue such optimisations further for the time being. Instead,  we focus on the optimal truncation problem in the remainder of the paper.

\subsection{Truncation based on the residual}\label{ss:residual}

We have observed that the coefficients of frame approximations alone do not offer sufficient information about optimal truncation or about the accuracy of the approximation. Fortunately, there is a simple alternative:  the size of the residual of the least squares problem. The residual is of little use when computing an interpolant using a basis. In that case, it is always zero \cite[Chapter XI]{davis1975interpolation}. If all goes well, interpolation yields a well-conditioned square linear system that can be solved to high accuracy, regardless of the accuracy of the corresponding function approximation (i.e., independent of the error between the interpolation points). In the setting of this paper, i.e., with oversampling, the residual does correspond more closely to the approximation error. Recall \eqref{eq:sampling_error} and note that $f - \mathcal T_N\z = f - \sum_{k=1}^N z_k \phi_k$ is the difference between $f$ and its approximation. The discrete residual vector is the sampled version of the continuous residual function, and it remains to verify the connection between $\Vert A \c - \b\Vert = \Vert f - \mathcal T_N \z \Vert_M$ and $\Vert f - \mathcal T_N \z \Vert_{\H}$.

A good approximation in the continuous norm  immediately implies small error in the $M$-norm as well, if the functionals $l_{m,M}$ are bounded. If we denote the sampling operator by
\begin{equation}\label{eq:sampling_operator}
 {\mathcal S}_M(f) : \G \to \mathbb{C}^M : f \to \{ l_{m,M}(f) \}_{m=1}^M,
\end{equation}
then indeed we have
\[
 \Vert f - \mathcal T_N \z \Vert_M = \Vert {\mathcal S}_M (f - \mathcal T_N \z) \Vert \leq S_M \Vert f - \mathcal T_N \z \Vert_{\H},
\]
where $S_M = \|{\mathcal S}_M\|$. Note that if \eqref{eq:richness} holds, we expect $S_M \leq B'$.

However, the question at hand is the converse. Does a small discrete residual imply a small continuous residual? That is, can we find a bound of the form:
\[
 \Vert f - \mathcal T_N \z \Vert_{\H} \leq C_M \Vert f - \mathcal T_N \z \Vert_M,
\]
with some constant $C_M < \infty$? The answer to this question is, unfortunately, \emph{no}, since for most spaces $\G$ and $\H$ one can find a function $f \in \G$ satisfying $\Vert f \Vert_{\H} > 0$ but such that ${\mathcal S}_M f = 0$. Yet, this is a problem for \emph{all} adaptive numerical methods based on a discretisation. Indeed, even for Chebyshev approximation, the function $f(x) = e^x + T_{100}(x)$ would produce exactly the same convergence plot as in the top-left panel of Figure~\ref{fig:1dexample}. Truncation based on coefficient size at $N < 101$ would incur $\mathcal O(1)$ error of the approximation, invisible to the discretisation.\footnote{The constructor in Chebfun samples the function at a few random points in order to avoid such scenarios with a certain probability. We will adopt a similar procedure shortly.}

Fortunately, since we oversample, we can at least make a statement about $C_M$ in the asymptotic limit $M \to \infty$ as a direct consequence of \eqref{eq:richness}.

\begin{lemma}\label{lem:norm_equivalence}
 Let $f\in G \subset \H$ and $\z \in \mathbb{C}^N$, let $\Phi_N \triangleq \{ \phi_k \}_{k=1}^N$ be the truncation of a frame $\Phi$ for $\H$ with all its elements satisfying $\phi_k \in \G$, and let a sampling operator ${\mathcal S}_M : \G \to \mathbb{C}^M$ be given by~\eqref{eq:sampling_operator}. If the associated norm $\Vert f \Vert_M \triangleq \Vert S_Mf \Vert$ satisfies \eqref{eq:richness}, then
\[
 \Vert f - \mathcal T_N \z \Vert_{\H} \leq \frac{1}{\sqrt{A'}} \sqrt{\liminf_{M \to \infty}\Vert f  - \mathcal T_N \z \Vert_M}.
 \]
\end{lemma}
\begin{proof}
 This is a consequence of \eqref{eq:richness} for the function $g \triangleq  f - \mathcal T_N \z$, noting that $g \in G$.
\end{proof}

We have already made one crucial assumption about $M$, namely that $M = \Theta^\epsilon(N,\theta)$ follows a stable sampling rate as a function of $N$. Loosely speaking, this condition guarantees that $\Vert f \Vert_M$ carries sufficient information to recover $f \in \H$ if the latter can be well approximated in $\H_N \triangleq \myspan \Phi_N$. Hence, we are within the asymptotic regime in which the discrete residual is a good substitute for the continuous residual. However, as long as $f$ is unresolved ($N$ is too small), we can not rely on the residual alone. We will have to augment any check of the residual with other checks to verify the onset of convergence.

\subsection{A stopping criterion}\label{ss:stopcriterion}

After solving the approximation problem $A \c^\epsilon=\b$, the computable tools at our disposal are the coefficient norm $\|\c^\epsilon\|$, the right hand side norm $\|\b\|$ and the residual norm $\|A \c^\epsilon - \b\|$. Recall that $\|\b\|= \| \mathcal S_M f \| = \|f\|_M$ and $\|A \c^\epsilon - \b\| = \| f - \mathcal T_N \c^\epsilon_N \|_M$.

We aim for a combination of \eqref{eq:criterion_part1} and \eqref{eq:criterion_part2}, e.g.,
\[
\|f-f_N\|_\H \leq \delta \|f\|_\H \quad \mbox{and} \quad  \| \c_N \| \leq \mu \|f\|_\H.
\]

\subsubsection{Coefficient norm versus residual norm}

We set out to illustrate that the coefficient norm and the residual norm convey similar information, yet the residual norm is preferable to use in practice.

\begin{figure}[tb]
	\centering
	\resizebox{.6\columnwidth}{!}{\input{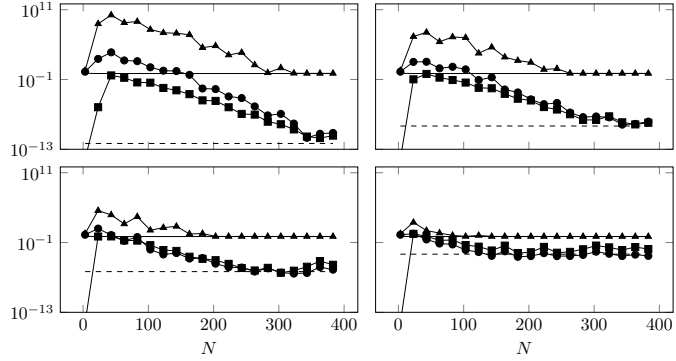}}
	\caption{Evolution of the $\H$-norm approximation error (dots), residual (squares) and coefficient size (triangles) of the Fourier extension approximation to $f(x)=e^{\cos(8\pi x)}$ for $\epsilon=10^{-12},10^{-9},10^{-6},10^{-3}$ (indicated with dashed black horizontal line). The asymptotic limit $\|f\|/\sqrt{A}$ for the coefficient size is indicated with the solid black horizontal line.}\label{fig:groupplotexpcos8pix}
\end{figure}

%\begin{figure}[tb]
%	\centering
%	\resizebox{.8\columnwidth}{!}{\input{img/groupplotfraction.tikz}}
%	\caption{\label{fig:groupplotfraction} Evolution of the $L^2$ approximation error (blue dots), residual (red squares) and coefficient size (brown crossed dots) of  the Fourier extension approximation to $f(x)=\frac{1}{x-1.1}$ for $\epsilon=10^{-12},10^{-9},10^{-6},10^{-3}$.}
%\end{figure}

It is clear from the previous subsections that it seems a good idea to consider a residual-based criterion
\begin{equation}\label{eq:residualbasedcriterion}
\|A \c^\epsilon - \b\| \leq \delta \| \b \|.
\end{equation}
Indeed, the residual is equivalent to the pointwise approximation error $\|f-\mathcal T_N\c\|_M$ and, once $M$ is sufficiently large, it tends to (a multiple of) $\|f-T_N\c\|_\H$ by Lemma~\ref{lem:norm_equivalence}. Analogously, $\|\b\|=\|f\|_{M}$ tends to $\|f\|_\H$. These expressions return in \eqref{eq:criterion_part1}. We can thus approximately replace the true stopping criterion \eqref{eq:criterion_part1} by the discrete and computable condition \eqref{eq:residualbasedcriterion}, at least for large $M$. For simplicity, we use the same constant $\delta$ in both inequalities.

As is evident from the coefficient bound \eqref{eq:coefficient_bound}, the coefficient norm can grow large in the regime before the onset of convergence. From \cite[Corollary 3.6]{adcock2017frames}, we also obtain (under the stable sampling rate) that
\[
 \limsup_{N\rightarrow\infty}\|\c^\epsilon\| \leq \frac{1}{\sqrt A} \|f\|,
\]
where $A$ is the lower frame bound of the frame.\footnote{Note that the coefficients can not be bounded for all $f \in \H$ if $\Phi_N$ is the truncation of an infinite set $\Phi$ that does not have a lower frame bound on $\H$. On the other hand, for any given function $f$, approximations with bounded coefficients may exist without any restrictions on $\Phi_N$ and, if so, they will also be found numerically. Unfortunately, in that case one may also find functions $f$ for which \eqref{eq:criterion_part2} can not be satisfied. In other words: the frame condition guarantees success for all $f \in H$, but absence of the frame condition does not prevent success for many functions $f \in \H$.}
The stability requirement  \eqref{eq:criterion_part2} can be meaningfully replaced by a corresponding condition involving computable quantities
\begin{equation}\label{eq:coefficientbasedcriterion}
	\| \c^\epsilon \| \leq \mu \| \b \|.
\end{equation}
Here, we again assume that $M$ is sufficiently large. 

Figure~\ref{fig:groupplotexpcos8pix}  illustrates the evolution of the coefficient norm (triangles) for increasing $N$ along with the  $\H$-norm error (dots) and the residual  (squares) for $\epsilon=10^{-12}$ (top left), $\epsilon=10^{-9}$ (top right), $\epsilon=10^{-6}$ (bottom left), and $\epsilon=10^{-3}$ (bottom right). The residual tends to the truncation parameter $\epsilon$ (dashed horizontal line) and the  coefficient norm to ${\|f\|_\H}/{\sqrt A} $ (full horizontal line). This is in full agreement with the stated theory.

The figure illustrates several important features. Firstly, the residual alone is no good estimate for the $\H$-norm error for small $N$ (and small $M$). The error of the approximation is invisible to the discretisation. Secondly, the coefficient norm does indeed attain the worst case $\mathcal O(\epsilon^{-1})$, suggesting that the second criterion \eqref{eq:coefficientbasedcriterion} might be needed. However, the figure also shows that the coefficient norm reaches $\|f\|/\sqrt A$ approximately for the same value of $N$ where the residual reaches $\epsilon$. Loosely speaking, if $\delta\approx\epsilon$, criteria \eqref{eq:residualbasedcriterion} and \eqref{eq:coefficientbasedcriterion}  convey similar information with a well chosen $\mu$. This is to be expected, as the coefficient bound given by \eqref{eq:coefficient_bound} indicates that the coefficient size is related to the residual norm. 

However, a suitable choice of $\mu$ depends on the frame dependent constant $A$, which may not be known in practice. In contrast, it is far easier to ensure that $A' \approx B' \approx 1$ by weighting the sample points as a Riemann sum, independently of what $\Phi_N$ is, to ensure close correspondence between the $\H$ and $M$ norms. Also, the coefficient norm tends to decrease very slowly towards its limit, creating a risk of obtaining an optimal approximation size that is way too large if $\mu$ is chosen just a bit too small.  We therefore prefer the residual-based criterion.

Another important piece of information conveyed by Figure~\ref{fig:groupplotexpcos8pix} is that choosing $\delta$ and $\epsilon$ far away leads to other problems. 
Firstly, if $\delta$ is too large with respect to $\epsilon$, and using the residual-based criterion \eqref{eq:residualbasedcriterion} only, we risk that the coefficient norm is high: $\|\c^\epsilon\|=\mathcal O(\delta/\epsilon)$. 
E.g., in Figure~\ref{fig:groupplotexpcos8pix} we see that the choice of $\epsilon=10^{-12}$, $\delta=10^{-3}$  and the residual-based criterion leads to the optimal value optimal $N=150$, but with a coefficient norm of $\mathcal O(10^9)$. 
Secondly, if $\epsilon$ is too large with respect to $\delta$, we can never solve the underlying linear system at the accuracy required by $\delta$. We conclude that $\epsilon$ and $\delta$ should be chosen comparable in size, and that $\epsilon$ should be small enough to allow solutions that reach the desired accuracy $\delta$ as $N$ increases. We will choose $\epsilon$ slightly smaller than $\delta$.

\subsubsection{Residual-based stopping criterion for all $N$}\label{sss:residualbased}

We adopt the residual-based criterion \eqref{eq:residualbasedcriterion} and discard the coefficient-based criterion \eqref{eq:coefficientbasedcriterion}. However, since \eqref{eq:residualbasedcriterion} is valid only in the regime of convergence, we need to augment the conditions to avoid returning an unacceptable solution for small $N$ that happens to pass the test by chance. Following the example of Chebfun, we do so by evaluating in a few additional random points.

\begin{figure}
	\centering
	\resizebox{.9\columnwidth}{!}{\input{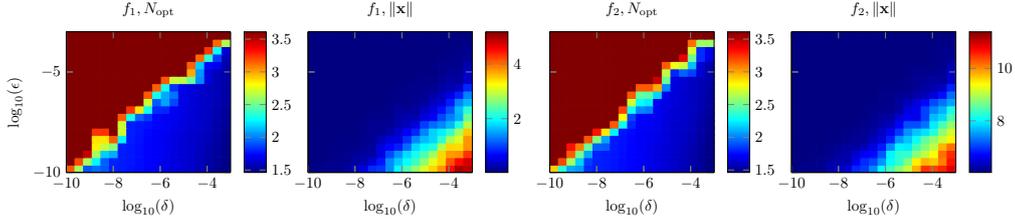}}
	\caption{Illustration of the optimal value of $N$ and the coefficient norm for varying $\delta$ and $\epsilon$ according to the stopping criterion of \S\ref{sss:residualbased}. Shown are the results for $f_1(x) = e^{\cos(8\pi)}$ (first two panels) and $f_2(x) = 10^6f_1(x)$ (last two panels). The functions are approximated on $[0,1/2]$ using Fourier series on $[0,1]$.  In both cases, the left panel depicts $N$ and the right panel shows $\Vert \x \Vert$, both in logarithmic scales. Even if no sufficient solution is found, the search for optimal $N$ stops at $N=4096$. The exact algorithm used is Algorithm~\ref{alg:bisection} with $\delta'=\delta$, $Q=3$.}\label{fig:ed}
\end{figure}

Given constants $\delta$, $\delta'$ and $\epsilon$, we accept a solution $\c^\epsilon_N$ if all of the following conditions are met:
\begin{enumerate}
 \item $\|A \c^\epsilon - \b\| \leq \delta \| \b \|$, and
 \item $|f(t_i) - f_N(t_i)| \leq \delta'\|\b\|$, $i=1,\ldots,Q$, where $\{ t_i \}_{i=1}^Q$ are random points in $\Omega$.
\end{enumerate}
The conditions are verified in order.

As argued in \S\ref{ss:residual}, we can discard a check on the coefficient norm if $\delta$ and $\epsilon$ are close. It is the second condition that, with some probability of success, catches the case where the residual is small in the $M$-norm, yet the approximation is large in between the sampling points. This is achieved by verifying for the desired accuracy in a few random points. For all examples in this paper, this condition is active only at very small values of $N$ and $M$, but it may also detect problems with degenerate cases such as $f(x)=e^x+\phi_{100}(x)$ where $\phi_{100}\notin\myspan\{\Phi_{99}\}$.

Some care is needed for the choice of $\delta'$. On the one hand, since 
%\begin{align}
%	\|f-f_N\|_{\H} \leq  \delta'|\Omega|
%\end{align}
the second restriction might be more stringent than the first if $\delta'$ is too small. On the other hand, if $\delta'$ is too large we might approximate $f(x)=e^x+\delta\phi_{100}(x)$ inaccurately. In the examples below, it is safe to choose $\delta=\delta'$.

These criteria are illustrated in Figure~\ref{fig:ed}. We compare the optimal $N$ and the coefficient norm $\|\c^\epsilon_N\|$ for $\epsilon,\delta\in[10^{-10},10^{-3}]$, with $\delta'=\delta$ and $Q=3$. If $\epsilon$ is too large compared to $\delta$, no optimal value of $N$ can be found. In the figure, these solutions have the maximal value $N=4096$, and they can be found in the top left triangles in the first and third panel. If $\delta$ is large and $\epsilon$ is small, the coefficient norm grows large. This is visible in the bottom right corners of panels 2 and 4.

Figure~\ref{fig:ed} also illustrates the scale-invariance of the relative stopping criterion \eqref{eq:residualbasedcriterion}. The optimal value of $N$ found for the functions $f_1(x) = e^{\cos(8\pi)}$ and $f_2(x) = 10^6f_1(x)$ are comparable while the coefficient norm is clearly $10^6$ times larger for the latter function.

\begin{figure}
	\centering
	\resizebox{.5\columnwidth}{!}{\input{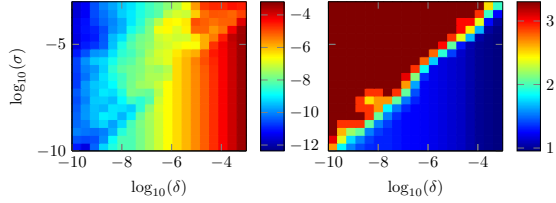}	}
	\caption{Approximation of a function with noise, $f(x) = e^x + \sigma\cos(2000\pi x)$, on $[0,1/2]$ using Fourier series on $[0,1]$. Shown are the residual in logarithmic scale, $\log_{10}(\|A_N\c^\epsilon_N-\b_N\|)$ (left panel), and the optimal value of $N$, $\log_{10}(N)$ (right panel), according to the criteria in \S\ref{sss:residualbased} and for varying $\delta$ and $\sigma$. The exact algorithm used is Algorithm~\ref{alg:bisection} with $\delta'=\delta$, $Q=3$ and $\epsilon=\delta/100$. The coefficient norm is large whenever $\delta$ is smaller than $\sigma$. Best results are obtained for $\delta\approx \sigma$, shown along the anti-diagonal.}\label{fig:sd}
\end{figure}

Finally, in Figure~\ref{fig:sd} we show results for a function with noise. We choose $f(x) = e^x + \sigma\cos(2000\pi x)$. In this case we need to choose $\delta$ appropriately, ideally comparable in size to the noise level. If $\delta$ is too small, then we try to approximate the noise, which results in large values of $N$. We do not include a description in this paper of methods for the detection of noise levels, instead relying on the user to choose $\delta$ appropriately.

\section{Algorithms}\label{s:algorithms}

We aim to find the smallest value of $N$ such that the approximation to a given function $f$ with $N$ degrees of freedom satisfies all conditions in \S\ref{sss:residualbased}:
\begin{equation}\label{eq:optimalN}
N(f;\delta,\delta',\epsilon,Q) = \argmin_{\substack{\|A \c^\epsilon_N - \b\| \leq \delta \| \b \|
\\ |f(t_i) - f_N(t_i)| \leq \delta'\|\b\|, i=1,\ldots,Q } } N.
\end{equation}
Recall that we assume a stable sampling rate $M=\Theta(N,\theta)$ and $\epsilon$ sufficiently small, but close to $\delta$. In all examples below, $M=2N$ was a sufficient rate.

\subsection{Brute force incremental strategy}

The simplest algorithm attempts each value of $N$ in order, until the optimal value $N(f; \delta, \delta',\epsilon,Q)$ is found. This approach yields the optimal value of $N$ in all cases. Obviously, it is very computationally expensive, as it takes the solution of $\OO(N_{\textrm{opt}})$ approximation problems to find the answer. Still, we have implemented this approach as a benchmark to compare with later on. The approach is specified in Algorithm~\ref{alg:incremental} where, as in Algorithm~\ref{alg:rapiddecay}, we use $\{A_i\}_{i=1}^\infty$ and $\{\b_i\}_{i=1}^\infty$ to denote the linear systems that need to be solved for a frame approximation using $i$ coefficients. The number of samples $M$ is implicit.

\begin{algorithm}
	\caption{The incremental strategy}\label{alg:incremental}
	{\bf Input:} Given $\delta,\delta', \epsilon,Q, \{A_i\}_{i=1}^\infty,\{\b_i=\mathcal S_if\}_{i=1}^\infty$ \\
	{\bf Output:} Optimal $N$
	\begin{algorithmic}
		\State $N = 1$
		\State Solve $A_N\c_N^\epsilon = \b_N$
		\While{ $\|A_N \c_N^\epsilon - \b_N\| >  \delta \| \b \|$ or  $|f(t_i) - f_N(t_i)| > \delta'\|\b_N\|$, $i=1,\ldots,Q$ }% \texttt{<some stopping criterion>}
		\State $N \gets N+1$
		\State Solve $A_N\c^\epsilon_N = \b_N$
		\EndWhile
		\State \Return $N$
	\end{algorithmic}
\end{algorithm}

Note that Algorithm~\ref{alg:incremental} is exactly like Algorithm~\ref{alg:rapiddecay}, but without the adaptive choice of weights. If the weighted approximation problem can be solved efficiently, then Algorithm~\ref{alg:rapiddecay} is not much more expensive than Algorithm~\ref{alg:incremental}.

\subsection{Bisection strategy}

The incremental strategy can be optimised with a relatively small modification that leads to substantially improved computational cost. We consider a method that consists of two stages. It is specified in Algorithm~\ref{alg:bisection}.

In the first stage, the method finds a size $\hat N$ that meets the stopping criterion, but which may not be optimal. This is achieved by doubling the number of degrees of freedom in every iteration until $\hat N$ is found. Hence, it takes $\OO(\log N_{\textrm{opt}})$ steps where $N_{\textrm{opt}} =N(f; \delta,\delta',\epsilon,Q)$ is the optimal value. In practice, we also stop at a maximal value of $N_{\textrm{max}}$ to ensure that the algorithm terminates.

In the second stage, the optimal step size is known to be between $N_h=\hat N$ and $N_l=\hat N/2$. The precise value is found using bisection. This takes an additional $\OO(\log N_{\textrm{opt}})$ steps. Hence, in total the scheme converges to an answer based on the solution of $\OO(\log N_{\textrm{opt}})$ approximation problems, rather than $N_{\textrm{opt}}$ approximation problems in Algorithm~\ref{alg:incremental}.

The first doubling stage is also employed in the software packages Chebfun \cite{driscoll2014fun,trefethen2013approximation} and ApproxFun~\cite{olver2018approxfun}. The second stage is different from the aforementioned packages, in which the coefficients can be truncated based on their decay and based on plateau detection to determine a noise level~\cite{Aurentz2017}. Plateau detection is impractical in our case, because we can only observe a plateau in the residuals for increasing number of degrees of freedom, and that is expensive to compute. In comparison, we have to solve a small number of additional approximation problems for values of $N$ between $N_h$ and $\frac{N_h}{2}$.

\begin{algorithm}
	\caption{Doubling and bisection}\label{alg:bisection}
	{\bf Input:} Given $\delta, \delta', \epsilon,Q, \{A_i\}_{i=1}^\infty,\{\b_i=\mathcal S_if\}_{i=1}^\infty$\\
	{\bf Output:} Optimal $N$
	\begin{algorithmic}
		\State $N = 1$
		\State Solve $A_N\c_N^\epsilon = \b_N$
		\While{ $\|A_N \c_N^\epsilon - \b_N\| >  \delta \| \b \|$ or  $|f(t_i) - f_N(t_i)| > \delta'\|\b_N\|$, $i=1,\ldots,Q$ }\Comment{Doubling}% \texttt{<some stopping criterion>}
		\State $N = 2N$
		\State Solve $A_N\c^\epsilon_N = \b_N$
		\EndWhile
		\State $\hat N\gets N$, $N_{h} \gets N$, $N_{l} \gets N/2$
		\While{$\left|N_h-N_l\right|>0$}\Comment{Bisection}
		\State $N = \left[\frac{N_l+N_h}{2}\right]$
		\State Solve $A_N\c_N^\epsilon = \b_N$
		\If{$\|A_N \c_N^\epsilon - \b_{N}\| >  \delta \| \b_{\hat N} \|$ or  $|f(t_i) - f_N(t_i)| > \delta'\|\b_{\hat N}\|$, $i=1,\ldots,Q$ }
		\State $N_l \gets N$
		\Else
		\State $N_h \gets N$
		\EndIf
		\EndWhile
		\State\Return $N$
	\end{algorithmic}
\end{algorithm}

In the description of the algorithm and the estimation of its cost, we have assumed for simplicity that $N$ can take all values between $1$ and $N_{\textrm{max}}$. In our implementation, we have merely assumed that there is a sequence of sets $\Phi_{j}$ with strictly increasing cardinality $N(j)$. This is considerably more general, but highly similar to the algorithm as stated here. The final computational cost does depend on the relation between $N(j)$ and $j$.

For example, in the case of a multivariate approximation with tensor product form, we may have that $N = \prod_{i=1}^d N_i$, using $N_i$ degrees of freedom in dimension $i$. If each $N_i$ is doubled in phase 1 of the algorithm, then $N$ increases by a factor of $2^d$. Alternatively, each $N_i$ can be multiplied by $2^{1/d}$, and subsequently rounded to the nearest integer. This results in roughly an increase by a factor of $2$ for the total number of degrees of freedom. Yet, even then, not all integers $N$ can be reached -- at least the prime numbers are always excluded. Care has to be taken in the implementation of the bisection algorithm to ensure that it terminates, and that only valid values of $N$ are used.

Another example is the approximation using a function dictionary of the form
\[
\Phi_N \triangleq \{ T_k(x) \}_{k=0}^{N_1-1} \cup \{ (1-x)^\alpha T_k(x) \}_{k=0}^{N_2-1}.
\]
This dictionary is well suited to approximate functions with an algebraic singularity of order $\alpha$ at $x=1$. Here, $N = N_1 + N_2$ is the sum of the lengths of the constituent sets (rather than their product as in the tensor product example). Depending on the application, it may be desirable to vary $N_1$ while keeping $N_2$ fixed, or one may vary $N_1$ and $N_2$ simultaneously. Our own implementation is configurable and allows both options, and in fact even allows the combination of tensor product dictionaries and weighted dictionaries, see Figure~\ref{fig:2dwFE}.

\subsection{On the monotonicity of convergence}

The bisection strategy is quite clearly more performant than the incremental strategy. However, there is an important qualitative difference: it is not a priori guaranteed that the bisection approach actually yields the optimal value of $N$. Assuming that random sampling is sufficient to detect all degenerate cases, the incremental approach always yields the minimal $N$, since it is based on the enumeration of all possibilities. The bisection approach may miss the optimal value, if the decrease of the error is not monotonic.

Say the optimal value in an approximation problem lies between $N_l = 16$ and $N_h=32$, and say, $N_{\textrm{opt}}=18$. If the bisection algorithm notices that convergence criteria are not satisfied for $N = (N_l+N_h)/2=24$, then it will continue to look for an optimal value in the interval $[25,32]$. Whatever it returns will be larger than $18$.

We explore the extent to which regularised frame-approximations such as those described in \S\ref{ss:frameapproximation} and \S\ref{ss:discrete_approximation} yield monotonic convergence rates. In the absence of regularisation, monotonicity is often guaranteed. Indeed, if the bases $\Phi_N$ form a nested sequence, i.e., $\myspan \Phi_{N-1} \subset \myspan \Phi_{N}$, then the best approximation error decays monotonically. This is the case if $\Phi_N$ arises from the truncation of an infinite set.

%For exact approximations, $\mathcal P_N$, the approximation spaces are nested in each other, i.e., $\H_{N-1} \subseteq \H_N$. From which follows monotone decrease of the approximation error:
\begin{lemma}
 Let $\{ \Phi_N \}_{N=1}^\infty$ form a nested sequence that is ultimately dense in a Hilbert space $\H$.  Let $f\in\H$, and let the orthogonal projection of $f$ on the truncated space $\H_N \triangleq \myspan \Phi_N$ be $\mathcal{P}_Nf$, then
	\begin{equation}
	\| f - \mathcal{P}_Nf \|_\H \leq \|f- \mathcal{P}_{N-1}f\|_\H.
	\end{equation}
\end{lemma}
\begin{proof}
The orthogonal projection in $\H$ onto $\H_N$ satisfies by construction
\[
 \|f-\mathcal{P}_Nf \|_\H \leq  \| f-\mathcal{T}_N\z \|_\H, \quad \forall \z\in\CC^N.
\]
Since $\H_{N-1}\subseteq \H_{N}$, this also holds for all $\mathcal T_{N-1}\z$, with $\z\in\CC^{N-1}$,  and more specifically for
\[
 \mathcal{P}_{N-1}f = \argmin_{\z\in\CC^{N-1}}\| f-\mathcal T_{N-1}\z \|_\H.
\]
\end{proof}

Monotonicity may be violated if the sequence is not nested, or if the discrete solution does not accurately reflect the best approximation. The latter may be a consequence of a coincidental bad choice of sampling points, but is affected more substantially by the regularisation.

Indeed, the lemma above does not hold for the regularised approximation spaces $\H^\epsilon_N \subset \H_N$ that result from the regularisation. Let \[ G_{M,N}=\{\ell_{k,M}(\phi_l)\}_{k,l=1}^{M,N}\]  be the discretisation matrix of a dictionary $\Phi_N$ using a sampling operator
\[
f\rightarrow \{\ell_{k,M}(f)\}_{k=1}^M,
\]
with singular value decomposition $G_{M,N} = U\Sigma V'$, where $U\in\CC^{M\times M}$, and $V\in\CC^{N\times N}$ are unitary, and $\Sigma\in\RR^{M\times N}$ is a diagonal matrix with positive entries $\sigma_k$. Using a truncated SVD to regularise the system (by discarding all singular values smaller than a given threshold $\epsilon$),
we define $\H^\epsilon_N$ as the span of the right singular vectors $v_k$ corresponding to singular values $\sigma_k$ larger than $\epsilon$: 
\begin{equation}\label{eq:defHepsilon}
\H^\epsilon_N\triangleq\myspan\{\mathcal{T}_N v_k|\sigma_k>\epsilon\}.
\end{equation}

Since $\H^\epsilon_{N-1}$ is not necessarily included in $\H^\epsilon_N$, even when $\H_{N-1} \subset \H_N$, a decrease in approximation error is not guaranteed. However, using the analysis in~\cite{Adcock2019}, we have the following:

\begin{theorem}\label{thm:decay}
	Let $f\in\H$, and the orthogonal projection of $f$ on $\H_N^\epsilon$ be $\mathcal{P}^\epsilon_Nf$, then
	\begin{equation}\label{eq:fna1_monotone}
	\|f-\mathcal{P}_N^\epsilon f\|_{\H} \leq \|f-\mathcal P^\epsilon_{N-1}f\|_{\H} + \sqrt\epsilon\|\c^\epsilon_{N-1}\|_{\ell^2(\CC^{N-1})},
	\end{equation}
	where $\c^\epsilon_{N-1}$ is the coefficient vector of $\mathcal P^\epsilon_{N-1}f$.
	\begin{proof}
		Equation \eqref{eq:fna1bound} says
		\begin{align*}
			\|f-\mathcal P^\epsilon_Nf\|_{\H} \leq \|f-\mathcal T_N\z_N\|_{\H} + \sqrt\epsilon\|\z_N\|_{\ell^2(\CC^N)},\quad \forall \z_N\in\CC^{N}.
		\end{align*}  
		Since this holds for any $\z_N\in\CC^N$, we can choose any $\z_N$ for which the last element is zero and obtain
		\begin{align*}
		\|f-\mathcal P^\epsilon_Nf\|_{\H} \leq \|f-\mathcal T_{N-1}\z_{N-1}\|_{\H} + \sqrt\epsilon\|\z_{N-1}\|_{\ell^2(\CC^{N-1})},\quad \forall \z_{N-1}\in\CC^{N-1},
		\end{align*}
		and more specifically
		\begin{align*}
		\|f-\mathcal P^\epsilon_Nf\|_{\H} \leq \|f-\mathcal P^\epsilon_{N-1}f\|_{\H} + \sqrt\epsilon\|\c\|_{\ell^2(\CC^{N-1})},
		\end{align*}
		where $\c$ is the coefficient vector of $\mathcal P^\epsilon_{N-1}f$.
	\end{proof}
\end{theorem}

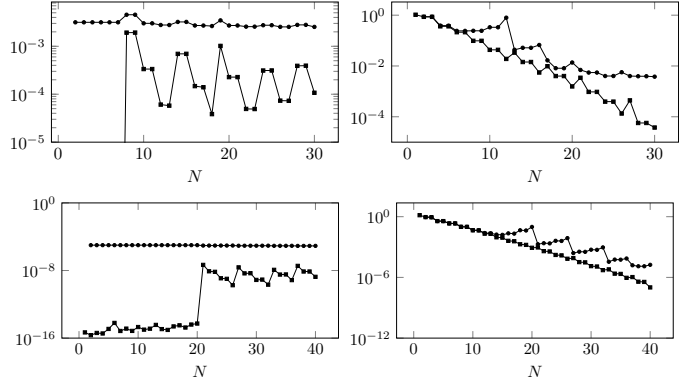
\begin{figure}[t]
	\centering
	\resizebox{.6\columnwidth}{!}{\begin{tikzpicture}
\begin{axis}[width={.5\textwidth}, height={.3\textwidth}, ymode={log}, xlabel={$N$}, cycle list name={mark list*}, ymin={1.0e-5}]
    \addplot+[mark size={1pt}]
        table[row sep={\\}]
        {
            \\
            1.0  nan  \\
            2.0  0.0031622776601686517  \\
            3.0  0.003162277660168407  \\
            4.0  0.0031622776601795354  \\
            5.0  0.0031622776601701466  \\
            6.0  0.0031622776601720323  \\
            7.0  0.0031622776601705526  \\
            8.0  0.004548488352751669  \\
            9.0  0.004550924221680684  \\
            10.0  0.0030113322357342145  \\
            11.0  0.003011038375253865  \\
            12.0  0.0027802439972758365  \\
            13.0  0.0027753017363339977  \\
            14.0  0.0032205786988655553  \\
            15.0  0.0032217275173666246  \\
            16.0  0.0027139538832712507  \\
            17.0  0.0027065661454614707  \\
            18.0  0.00264099485136476  \\
            19.0  0.003466991158699131  \\
            20.0  0.002714136701340927  \\
            21.0  0.0027143109313606936  \\
            22.0  0.002568607508633654  \\
            23.0  0.0025678988224356737  \\
            24.0  0.002748462301389981  \\
            25.0  0.0027494187070636706  \\
            26.0  0.0025396327664077953  \\
            27.0  0.002538371302561469  \\
            28.0  0.0027932427767196353  \\
            29.0  0.002796519029256276  \\
            30.0  0.0025356735996447588  \\
        }
        ;
    \addplot+[mark size={1pt}]
        table[row sep={\\}]
        {
            \\
            1.0  0.0  \\
            2.0  2.7111314086419895e-16  \\
            3.0  2.7755575615628914e-17  \\
            4.0  1.0927201054607919e-14  \\
            5.0  1.732224544821262e-15  \\
            6.0  4.056529615161281e-15  \\
            7.0  2.3951407705975006e-15  \\
            8.0  0.001925718130299314  \\
            9.0  0.0019292472936443777  \\
            10.0  0.0003349109825113192  \\
            11.0  0.00033480341450156236  \\
            12.0  6.07739337698435e-5  \\
            13.0  5.779948711928886e-5  \\
            14.0  0.0006946226278173582  \\
            15.0  0.0006960406785617306  \\
            16.0  0.0001475682600324329  \\
            17.0  0.00013964987454192997  \\
            18.0  3.847326702469705e-5  \\
            19.0  0.001017300166614476  \\
            20.0  0.00022741784157948578  \\
            21.0  0.0002275905409052205  \\
            22.0  4.9360498295154516e-5  \\
            23.0  4.909391953891741e-5  \\
            24.0  0.000311943251105286  \\
            25.0  0.0003128694419497089  \\
            26.0  7.34958474899464e-5  \\
            27.0  7.287416939007698e-5  \\
            28.0  0.00039098140707422046  \\
            29.0  0.0003941010187164424  \\
            30.0  0.00010729255441101336  \\
        }
        ;
\end{axis}
\end{tikzpicture}\begin{tikzpicture}
\begin{axis}[width={.5\textwidth}, height={.3\textwidth}, ymode={log}, xlabel={$N$}, cycle list name={mark list*}, ymin={1.0e-5}]
    \addplot+[mark size={1pt}]
        table[row sep={\\}]
        {
            \\
            1.0  nan  \\
            2.0  0.8654843864134045  \\
            3.0  0.8642521869924261  \\
            4.0  0.3886650759261323  \\
            5.0  0.3873746347441684  \\
            6.0  0.24097189174341474  \\
            7.0  0.24032421937853318  \\
            8.0  0.2429307800949695  \\
            9.0  0.24273266536853502  \\
            10.0  0.33367189182554335  \\
            11.0  0.33358800769651614  \\
            12.0  0.796104375697801  \\
            13.0  0.04335011793562799  \\
            14.0  0.051938698123972345  \\
            15.0  0.051899263496157935  \\
            16.0  0.06699381156710287  \\
            17.0  0.016610826259292694  \\
            18.0  0.00814700434345691  \\
            19.0  0.00812733044635059  \\
            20.0  0.013678098788709746  \\
            21.0  0.0069717467550417175  \\
            22.0  0.005507534407856153  \\
            23.0  0.005505825435417049  \\
            24.0  0.004023867518026303  \\
            25.0  0.004023289282300524  \\
            26.0  0.0056635505131323715  \\
            27.0  0.004003588068626097  \\
            28.0  0.003936487731649629  \\
            29.0  0.003936398043714622  \\
            30.0  0.0037912216688115894  \\
        }
        ;
    \addplot+[mark size={1pt}]
        table[row sep={\\}]
        {
            \\
            1.0  1.041777924567016  \\
            2.0  0.8611875518315093  \\
            3.0  0.8611875518315093  \\
            4.0  0.3639922923570177  \\
            5.0  0.3639922923570177  \\
            6.0  0.21378814361046994  \\
            7.0  0.21378814361046988  \\
            8.0  0.09780743050690999  \\
            9.0  0.0978074305069099  \\
            10.0  0.04316791299438972  \\
            11.0  0.043167912994390614  \\
            12.0  0.018972116282754012  \\
            13.0  0.03363176673674391  \\
            14.0  0.014322429289577616  \\
            15.0  0.014322410688960499  \\
            16.0  0.005534728378435206  \\
            17.0  0.009770813273807615  \\
            18.0  0.004006483669404381  \\
            19.0  0.004006431612070231  \\
            20.0  0.001565253393161617  \\
            21.0  0.0034380497903977723  \\
            22.0  0.0009540838838422082  \\
            23.0  0.0009540268597379237  \\
            24.0  0.000394537750268394  \\
            25.0  0.00039453486892699566  \\
            26.0  0.00013426796998499817  \\
            27.0  0.00044137099731062483  \\
            28.0  5.698947551005731e-5  \\
            29.0  5.6987371791506465e-5  \\
            30.0  3.7650788366393076e-5  \\
        }
        ;
\end{axis}
\end{tikzpicture}}
	\\
	\resizebox{.6\columnwidth}{!}{
	\begin{tikzpicture}
\begin{axis}[width={.5\textwidth}, height={.3\textwidth}, ymode={log}, xlabel={$N$}, cycle list name={mark list*}, ymax={1}, ymin={1.0e-16}]
    \addplot+[mark size={1pt}]
        table[row sep={\\}]
        {
            \\
            1.0  nan  \\
            2.0  1.000000000022194e-5  \\
            3.0  1.0000000000408247e-5  \\
            4.0  1.0000000000352998e-5  \\
            5.0  1.0000000001220857e-5  \\
            6.0  1.000000000627882e-5  \\
            7.0  1.0000000000698244e-5  \\
            8.0  1.0000000001328492e-5  \\
            9.0  1.0000000000712042e-5  \\
            10.0  1.000000000204548e-5  \\
            11.0  1.0000000001007718e-5  \\
            12.0  1.0000000001310295e-5  \\
            13.0  1.0000000003651851e-5  \\
            14.0  1.0000000001168727e-5  \\
            15.0  1.0000000000916602e-5  \\
            16.0  1.0000000002343753e-5  \\
            17.0  1.0000000003288931e-5  \\
            18.0  1.0000000002711777e-5  \\
            19.0  1.0000000004641421e-5  \\
            20.0  9.999999997664716e-6  \\
            21.0  8.947376451839537e-6  \\
            22.0  8.958898872847892e-6  \\
            23.0  8.958781734322799e-6  \\
            24.0  8.996459888421993e-6  \\
            25.0  8.997173235192532e-6  \\
            26.0  9.034766664833088e-6  \\
            27.0  8.545284921435665e-6  \\
            28.0  8.579896220453599e-6  \\
            29.0  8.581009370751421e-6  \\
            30.0  8.625181467175186e-6  \\
            31.0  8.626105048362948e-6  \\
            32.0  8.669086581228753e-6  \\
            33.0  8.32614031860056e-6  \\
            34.0  8.368245496531366e-6  \\
            35.0  8.369133526053145e-6  \\
            36.0  8.412969753329107e-6  \\
            37.0  8.165397102844816e-6  \\
            38.0  8.188247503989063e-6  \\
            39.0  8.188944660914614e-6  \\
            40.0  8.229923231440935e-6  \\
        }
        ;
    \addplot+[mark size={1pt}]
        table[row sep={\\}]
        {
            \\
            1.0  4.710277376051325e-16  \\
            2.0  2.219390173021042e-16  \\
            3.0  4.0824042601638266e-16  \\
            4.0  3.5300202911220376e-16  \\
            5.0  1.2208667421104631e-15  \\
            6.0  6.278840389516098e-15  \\
            7.0  6.982385683850684e-16  \\
            8.0  1.3283678574579974e-15  \\
            9.0  7.121551989798691e-16  \\
            10.0  2.0451368683336884e-15  \\
            11.0  1.0078681203721898e-15  \\
            12.0  1.3065440648386355e-15  \\
            13.0  3.652412509193777e-15  \\
            14.0  1.1670237262350449e-15  \\
            15.0  9.298442803701368e-16  \\
            16.0  2.4135866039163304e-15  \\
            17.0  3.2529696648268633e-15  \\
            18.0  1.789082409022188e-15  \\
            19.0  4.024996005892633e-15  \\
            20.0  5.294357468535218e-15  \\
            21.0  4.623783275375652e-8  \\
            22.0  8.080458341476266e-9  \\
            23.0  6.995156027097234e-9  \\
            24.0  1.1994572454549917e-9  \\
            25.0  1.0394892829289502e-9  \\
            26.0  1.81850147387554e-10  \\
            27.0  2.354302244498009e-8  \\
            28.0  4.574644095769507e-9  \\
            29.0  4.774450773617904e-9  \\
            30.0  7.70693948676522e-10  \\
            31.0  8.414311467295848e-10  \\
            32.0  2.1201440252639264e-10  \\
            33.0  1.20132636174313e-8  \\
            34.0  3.319897655248103e-9  \\
            35.0  3.419044078405876e-9  \\
            36.0  7.303612292212617e-10  \\
            37.0  3.5385642297521304e-8  \\
            38.0  7.623810487205792e-9  \\
            39.0  7.620655329521388e-9  \\
            40.0  1.7929281791706163e-9  \\
        }
        ;
\end{axis}
\end{tikzpicture}\begin{tikzpicture}
\begin{axis}[width={.5\textwidth}, height={.3\textwidth}, ymode={log}, xlabel={$N$}, cycle list name={mark list*}, ymin={1.0e-12}]
    \addplot+[mark size={1pt}]
        table[row sep={\\}]
        {
            \\
            1.0  nan  \\
            2.0  0.9092049192127587  \\
            3.0  0.8762720220822244  \\
            4.0  0.36623120772137757  \\
            5.0  0.36624824525423666  \\
            6.0  0.21578016121417776  \\
            7.0  0.21504805554283504  \\
            8.0  0.10033975560702996  \\
            9.0  0.09913930884817153  \\
            10.0  0.04621583025706167  \\
            11.0  0.045110210424931636  \\
            12.0  0.024196434502457335  \\
            13.0  0.023060308566558038  \\
            14.0  0.017114506649012068  \\
            15.0  0.01605262724687135  \\
            16.0  0.022592759903267787  \\
            17.0  0.02116178926900074  \\
            18.0  0.04575910314981472  \\
            19.0  0.04318452207180578  \\
            20.0  0.09737561373537611  \\
            21.0  0.001912337252908515  \\
            22.0  0.002389372000357006  \\
            23.0  0.002289852508692489  \\
            24.0  0.0041378330530184426  \\
            25.0  0.003977138602174687  \\
            26.0  0.007636186673226498  \\
            27.0  0.00025457041495471476  \\
            28.0  0.0003458654578362474  \\
            29.0  0.00033551456550453364  \\
            30.0  0.000545365403716662  \\
            31.0  0.0005292047857989793  \\
            32.0  0.0009213609658851558  \\
            33.0  3.504050871156097e-5  \\
            34.0  5.759871305685934e-5  \\
            35.0  5.6410102490701476e-5  \\
            36.0  7.358994154328178e-5  \\
            37.0  1.5744892819262427e-5  \\
            38.0  1.2613956570864258e-5  \\
            39.0  1.2595942044838137e-5  \\
            40.0  1.7501986278226516e-5  \\
        }
        ;
    \addplot+[mark size={1pt}]
        table[row sep={\\}]
        {
            \\
            1.0  1.4006410444073225  \\
            2.0  0.9091961206559691  \\
            3.0  0.8762636381277071  \\
            4.0  0.3661555943856442  \\
            5.0  0.36617729570804874  \\
            6.0  0.2156859729592631  \\
            7.0  0.21496419717641968  \\
            8.0  0.09982981350432427  \\
            9.0  0.09865411595789772  \\
            10.0  0.045109609230014974  \\
            11.0  0.044067631356666734  \\
            12.0  0.021055484243983267  \\
            13.0  0.020096969956032086  \\
            14.0  0.009129410878016749  \\
            15.0  0.008491835659703658  \\
            16.0  0.004125762837182467  \\
            17.0  0.0037157201391560374  \\
            18.0  0.001872662273331427  \\
            19.0  0.001634616335606987  \\
            20.0  0.0008477011560290287  \\
            21.0  0.0009123046343704618  \\
            22.0  0.0003921318933063993  \\
            23.0  0.00036392932026648237  \\
            24.0  0.0001654977335282941  \\
            25.0  0.00014953473726747027  \\
            26.0  6.991850211787808e-5  \\
            27.0  8.115599622523495e-5  \\
            28.0  3.2483973315068916e-5  \\
            29.0  3.0422002507023087e-5  \\
            30.0  1.2984518691062404e-5  \\
            31.0  1.1890934175778508e-5  \\
            32.0  5.215865440357441e-6  \\
            33.0  6.19790615338821e-6  \\
            34.0  2.35221607023676e-6  \\
            35.0  2.2080475172716636e-6  \\
            36.0  9.074762070247921e-7  \\
            37.0  1.076366668827648e-6  \\
            38.0  3.7682755181918756e-7  \\
            39.0  3.563252622526324e-7  \\
            40.0  1.0224194720897984e-7  \\
        }
        ;
\end{axis}
\end{tikzpicture}
	}
	\caption{\label{fig:notmonotone} Illustration of Theorems~\ref{thm:decay} and~\ref{thm:decay2} using a Fourier extension frame approximation for $f(x)=1$ (left panels) and $f(x) = \exp(\cos(5x))$ (right panels), with truncation parameter $\epsilon=10^{-5}$. Top row: approximation error $\|f-\mathcal{P}_N^\epsilon f\|_{\H}$ (squares) and the upper bound \eqref{eq:fna1_monotone} (dots). Bottom row: Approximation error $\|f-\mathcal{P}_{M,N}^\epsilon f\|_{\H}$ (squares) and the upper bound \eqref{eq:fna2_monotone}, here simplified to be $\|f-\mathcal{P}_{M',N-1}^\epsilon f\|_{\H}+\epsilon\|c\|_{\ell^2(\CC^{N-1})}$ (dots). Convergence is largely monotonic, but not completely, and no convergence occurs below the regularization threshold.}
\end{figure}

Therefore, bisection can not guarantee the most optimal solution, but it will give one with an error that is close to the optimal error.

The previous theorem is based on the bound \eqref{eq:fna1bound} using orthogonal projections. A similar result can be stated for the second bound \eqref{eq:fna2bound}, which applies to the discrete setting using samples.

\begin{theorem}\label{thm:decay2}
	Assume that $M$ and $M'$ satisfy the stable sampling rate for $N$ and $N-1$ respectively, i.e., $M\geq\Theta^\epsilon(N,\theta)$ and $M'\geq\Theta^\epsilon(N-1,\theta)$. Denote by $\mathcal{P}_{M,N}^\epsilon f$ the $\epsilon$-regularised solution of the discrete least squares problem, and let $A'$ be the lower bound constant in \eqref{eq:richness}, then
	\begin{equation}\label{eq:fna2_monotone}
	\|f-\mathcal{P}_{M,N}^\epsilon f\|_{\H} \leq \|f-\mathcal P^\epsilon_{M',N-1}f\|_{\H} +\frac{\theta}{A'}\left(
	 \|f-\mathcal \mathcal P^\epsilon_{M',N-1}f\|_M+\epsilon\|\c^\epsilon_{N-1}\|_{\ell^2(\CC^{N-1})}\right),
	\end{equation}
	where $\c^\epsilon_{N-1}$ is the coefficient vector of $\mathcal P^\epsilon_{M',N-1}f$.
\begin{proof}
The stable sampling rate is satisfied, so by \eqref{eq:bound}
\begin{align*}
\|f-\mathcal P^\epsilon_{M,N}f\|_{\H}\leq \|f-\mathcal T_N \z_N\|_{\H}+\frac{\theta}{A'}(\|f-\mathcal T_N \z_N\|_M+\epsilon\|\z_N\|), \quad\forall\z_N\in\CC^N. 
\end{align*}
Since this holds for any $\z_N\in\CC^N$, we can choose a vector $\z_{N-1} \in \CC^{N-1}$ and append a zero element at the end. This leads to
\begin{align*}
\|f-\mathcal P^\epsilon_{M,N}f\|_{\H}\leq \|f-\mathcal T_{N-1} \z_{N-1}\|_{\H}+\frac{\theta}{A'}(\|f-\mathcal T_{N-1} \z_{N-1}\|_M+\epsilon\|\z_{N-1}\|), \quad\forall\z_{N-1}\in\CC^{N-1}. 
\end{align*}
We obtain the result by choosing $\z_{N-1}$ equal to $\c^\epsilon_{N-1}$, the coefficient vector of $P^\epsilon_{M',N-1}f$.
\end{proof}
\end{theorem}
The theorem shows that the sequence is almost monotonic. Lack of monotonicity can be caused by the two rightmost terms in \eqref{eq:fna2_monotone}. The first of those, $\|f-\mathcal \mathcal P^\epsilon_{M',N-1}f\|_M$, is expected to be on the same order as $\|f-\mathcal P^\epsilon_{M',N-1}f\|_{\H}$ in the regime of convergence, i.e., for sufficiently large $N$. The second term may cause minor jumps in the approximation error on the order of $\epsilon$.

The behaviour of the $\H$-norm error and the upper bounds of Theorems~\ref{thm:decay} and~\ref{thm:decay2} are illustrated in Figure~\ref{fig:notmonotone}. Theorem~\ref{thm:decay} is illustrated in the top row for the constant function $f(x)=1$ (left panel) and a generic smooth function (right panel), using $\epsilon = 10^{-5}$. Theorem~\ref{thm:decay2} is illustrated in the bottom row for the same two functions. In all cases, the function is approximated using a Fourier extension frame.

Consider the constant function first (left panels). The regularisation yields a sequence of function spaces $\H^\epsilon_1,\H^\epsilon_2,\dots$. Only the first of these actually contains constant functions, the subsequent spaces do not. Hence, the error for $f(x)=1$ is machine precision at $N=1$, but is significantly larger afterwards. The upper bounds of the theorems are on the order of $\sqrt{\epsilon}$ and $\epsilon$ respectively. The true approximation error stays below this bound, but as expected no convergence can be achieved beyond the accuracy determined by the regularisation parameter.

The results for the smooth function show nearly monotonic behaviour, though not exactly, both for the upper bounds and the true residuals. The bounds allow for jumps on the order of $\sqrt{\epsilon}$ (top) and $\epsilon$ (bottom). For Theorem~\ref{thm:decay2}, we have simplified the bound to $\|f-\mathcal{P}_{M',N-1}^\epsilon f\|_{\H}+\epsilon\|\c\|_{\ell^2(\CC^{N-1})}$. Jumps are still present, but less visible and in the range allowed by the bound.

\subsection{Numerical experiments}

\subsubsection{A univariate Fourier extension example}

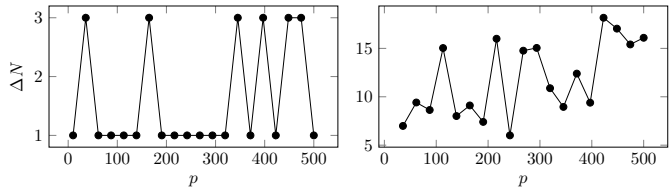
\begin{figure}[tb]
	\centering
	\resizebox{.6\columnwidth}{!}{
		\begin{tikzpicture}
\begin{axis}[cycle list name={mark list*}, xlabel={$p$}, width={.5\textwidth}, height={.3\textwidth}, ylabel={$\Delta N$}]
    \addplot+
        table[row sep={\\}]
        {
            \\
            10.0  1.0  \\
            35.78947368421053  3.0  \\
            61.57894736842105  1.0  \\
            87.36842105263158  1.0  \\
            113.1578947368421  1.0  \\
            138.94736842105263  1.0  \\
            164.73684210526315  3.0  \\
            190.52631578947367  1.0  \\
            216.3157894736842  1.0  \\
            242.10526315789474  1.0  \\
            267.89473684210526  1.0  \\
            293.6842105263158  1.0  \\
            319.4736842105263  1.0  \\
            345.2631578947368  3.0  \\
            371.05263157894734  1.0  \\
            396.8421052631579  3.0  \\
            422.6315789473684  1.0  \\
            448.4210526315789  3.0  \\
            474.2105263157895  3.0  \\
            500.0  1.0  \\
        }
        ;
\end{axis}
\end{tikzpicture}
		\begin{tikzpicture}
\begin{axis}[xlabel={$p$}, width={.5\textwidth}, height={.3\textwidth}, cycle list name={mark list*}]
    \addplot
        table[row sep={\\}]
        {
            \\
            35.78947368421053  6.9870532457326835  \\
            61.578947368421055  9.409257412447554  \\
            87.36842105263158  8.634024111792124  \\
            113.15789473684211  15.036657040075356  \\
            138.94736842105263  8.0098128045355  \\
            164.73684210526315  9.104901849224053  \\
            190.52631578947367  7.403545928681794  \\
            216.31578947368422  15.997153427436588  \\
            242.1052631578947  6.0107794143543805  \\
            267.89473684210526  14.762212375781694  \\
            293.6842105263158  15.048270120864638  \\
            319.47368421052636  10.885589071478755  \\
            345.2631578947368  8.952026891458086  \\
            371.05263157894734  12.393216873254435  \\
            396.8421052631579  9.389097561827748  \\
            422.63157894736844  18.16466643989254  \\
            448.4210526315789  17.039985285225004  \\
            474.2105263157895  15.406638930167317  \\
            500.0  16.09555578680503  \\
        }
        ;
\end{axis}
\end{tikzpicture}
	}
	\caption{Adaptive approximation of $f(x) = \cos(px)$ on $[-1,1]$ using Fourier series on $[-2,2]$ with $\delta=\epsilon=10^{-10}$ and $Q=3$. Left: The difference $N_{\textrm{bis}}-N_{\textrm{inc}}$ between the optimal value for $N$ obtained by the bisection algorithm compared to the brute-force incremental algorithm, as a function of $p$. The difference is very small. Right: the ratio of the computation time of the adaptive bisection approach over the time it takes to approximate $f$ once using the optimal number of degrees of freedom. To avoid the influence of the random point evaluation, the median is taken over 7 experiments for every data point.}\label{fig:adaptiveOptimalN}
\end{figure}

We start with a simple univariate example. In Figure~\ref{fig:adaptiveOptimalN},  we approximate $f(x) = \cos(px)$ on $[-1,1]$ using Fourier series on $[2,2]$, i.e., using Fourier extension. Larger values of $p$ render the function more oscillatory, and hence should lead to larger optimal values of $N$. We choose $p\in[0,500]$, $\delta=10^{-10}$, $\epsilon=10^{-12}$, $Q=3$.

We compare the optimal $N$ given by the incremental approach and the bisection approach. The bisection approach results in a slight overestimation of the true optimal $N$, obtained by the greedy approach. As seen in the figure, the difference is in this example never larger than $3$ degrees of freedom. For $p=500$, the optimal value of $N$ is $642$.

We also compare the computational cost of the bisection approach with that of a single function approximation using the optimal number of degrees of freedom. Let $t_{biss}$ be the time it takes to find $N_{\textrm{opt}}$, and let $t_{N_{\textrm{opt}}}$ be the time to compute an approximation with $N_{\textrm{opt}}$ degrees of freedom. The figure shows the ratio
\[
 \frac{t_{\textrm{biss}}}{t_{N_{\textrm{opt}}}}.
\]
We expect this ratio to be bounded by $\mathcal O(\log t_{N_{\textrm{opt}}})$ since the bisection approach computes only logarithmically many additional function approximations. In this example, the computational cost is cubic in $N$ for each approximation problem. The approximations with a large $N$ therefore take a much larger fraction of the total clock time than those with small $N$. The implication is that the logarithmic growth in the number of iterations is not actually visible in the figure.

\subsection{Two-dimensional spectral approximation of a singular function on a non-rectangular domain}

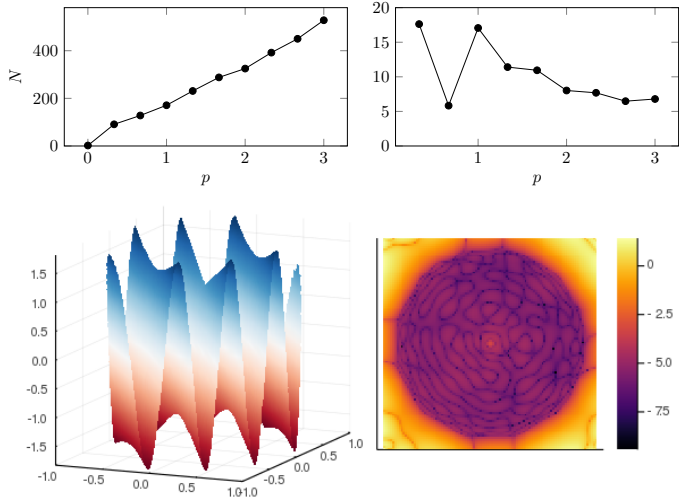
\begin{figure}[tb]
 	\centering
 	\resizebox{.6\columnwidth}{!}{\begin{tikzpicture}
\begin{groupplot}[group style={group size={2 by 1}}, cycle list name={mark list*}]
    \nextgroupplot[width={.5\textwidth}, height={.3\textwidth}, xlabel={$p$}, ylabel={$N$}, ymin={0}]
    \addplot
        table[row sep={\\}]
        {
            \\
            0.0  2.0  \\
            0.3333333333333333  91.0  \\
            0.6666666666666666  128.0  \\
            1.0  171.0  \\
            1.3333333333333333  231.0  \\
            1.6666666666666667  288.0  \\
            2.0  325.0  \\
            2.3333333333333335  392.0  \\
            2.6666666666666665  450.0  \\
            3.0  528.0  \\
        }
        ;
    \nextgroupplot[legend cell align={left}, width={.5\textwidth}, height={.3\textwidth}, ymin={0}, ymax={20}, xlabel={$p$}]
    \addplot
        table[row sep={\\}]
        {
            \\
            0.3333333333333333  17.624151941702838  \\
            0.6666666666666666  5.830572089075409  \\
            1.0  17.062438474565024  \\
            1.3333333333333333  11.406355335230637  \\
            1.6666666666666667  10.946076115938322  \\
            2.0  8.02175879708592  \\
            2.3333333333333335  7.689665891104797  \\
            2.6666666666666665  6.4786990200836465  \\
            3.0  6.792671335948985  \\
        }
        ;
\end{groupplot}
\end{tikzpicture}}
 	\includegraphics[width=.6\columnwidth]{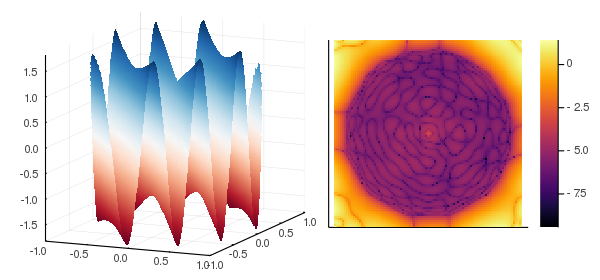}
 	\caption{We approximate the singular function \eqref{eq:example_2d} on the disk with centre $[0,0]$ and radius $0.9$ using the weighted Fourier Extension frame \eqref{eq:2dframe}. Upper left: The optimal number of degrees of freedom found by the adaptive approach using $\delta=10^{-6}$ and $\epsilon=10^{-8}$. 
 	Upper right: the ratio of the timings of the adaptive approach w.r.t. the timings of an approximation using the same number degrees of freedom as the optimal solution.. Lower left and right: the approximation and the $\log_{10}$ of the uniform error of the approximant found with $p=3,\delta=10^{-6},\epsilon=10^{-8}$.}\label{fig:2dwFE}
\end{figure}

We conclude with a more involved example that combines several difficulties. We choose a function that has an algebraic point singularity in 2-D,
\begin{equation}\label{eq:example_2d}
 f(x,y) = \cos(p\pi(x+y)) + \sqrt{x^2+y^2}\sin(1+p\pi(x+y)).
\end{equation}
As before, the parameter $p$ controls the oscillatory nature of the function. Furthermore, we set out to approximate this function on a non-rectangular domain, a disk with radius $0.9$. To that end, first we choose a Fourier extension frame by restricting a Fourier series $\Psi \triangleq \{ \psi_{k,l} \}_{k,l \in \ZZ}$ on $[-1,1]^2$ to the disk. Next, we use that frame for the smooth and singular parts of the function:
\begin{equation}\label{eq:2dframe}
 \Phi \triangleq \Psi \cup \sqrt{x^2+y^2} \Psi = \{ \psi_{k,l}(x,y) \}_{k,l \in \ZZ} \cup \{\sqrt{x^2+y^2} \psi_{k,l}(x,y) \}_{k,l \in \ZZ}.
\end{equation}
We define a truncation $\Phi_N$ with increasing length $N$ in terms of the truncated frames $\Psi_{N_1}$ and  $\sqrt{x^2+y^2} \Psi_{N_2}$, with $N=N_1+N_2$. When increasing $N$, we alternate between increasing $N_1$ and $N_2$, i.e., $N_1 = 1+ \lceil \frac{N-1}{2} \rceil$ and $N_2 = 1 + \lfloor \frac{N-1}{2} \rfloor$. In turn, $\Psi_{N_i}$ has size $N_i^{(1)}\times N_i^{(2)}$ and we similarly alternate between increasing $N_i^{(1)}$ and $N_i^{(2)}$ when increasing $N_i$. The first set $\Psi_1$ has size $1 \times 1$.

In Figure~\ref{fig:2dwFE}, we have repeated the experiment of Figure~\ref{fig:adaptiveOptimalN}. As before, we increase $p$ to obtain adaptive approximations with higher $N$ and compare timings and optimal $N$. The example is completed with a plot of the optimal approximant (bottom left) and the $\log_{10}$ of the uniform error (bottom right) for $p=3$, $\epsilon=10^{-8}$, $\delta=10^{-6}$, $Q=3$.

%\section{Conclusion}\label{s:conclusion}
%
%\BLUE{We presented the main difference between frame and basis approximations: coefficient decay is not guaranteed in the former. If the decay rate is known, one can obtain the coefficients at this decay rate solving a single weighted least squares. If this is not the case, the decay can be obtained using Algorithm~\ref{alg:rapiddecay}. }
%
%\BLUE{However, to determine the optimal $N$, i.e., the smallest $N$ that corresponds with an accurate and stable frame approximation \eqref{eq:criterion_part1}-\eqref{eq:criterion_part2}, we don't need coefficient decay. To solve frame approximation problems accurately, oversampling is key. Therefore, we can use the residual al an indicator for convergence. The  stability indicator is the coefficient norm which leads to the stopping criterion in \S\ref{sss:residualbased}. The simplest approach to arrive at the optimal $N$ is iterating over all $N$ until the optimal value is obtained is not practical because of the cost. We present in Algorithm~\ref{alg:bisection} a more efficient algorithm,  that may not generate the optimal $N$ but will be close, see theorems~\ref{thm:decay},~\ref{thm:decay2}. }
%
%\BLUE{The algorithms presented here, we assumed all bases can be indexed with all integer values of $N$. Our implementation is considerably more general, see the example in Figure~\ref{fig:2dwFE}.}

\bibliographystyle{abbrv}%IMANUM-BIB
\bibliography{adaptive_library}

\end{document}